\documentclass{article}

\PassOptionsToPackage{%
	style	= alphabetic%
}{biblatex}

\usepackage[doi=false,isbn=false,url=false,
	hyperref=auto,
	sorting=nyt,
	maxnames=10,
	maxcitenames=3,
	backend=biber,
	texencoding=auto,
	giveninits=true,
	block=space,
]{biblatex}

\DefineBibliographyStrings{english}{%
	andothers = {et al}
} 
 
\DeclareFieldFormat{eprint:arxiv}{%
	\href{https://arxiv.org/abs/#1}{\emph{arXiv:\allowbreak #1}}}

\DeclareFieldFormat{eprint:mrnumber}{%
	\href{http://www.ams.org/mathscinet-getitem?mr=MR#1}{MR\allowbreak #1}}

\DeclareFieldFormat{eprint:jstor}{%
	\href{http://www.jstor.org/stable/#1}{JSTOR\allowbreak #1}}

\DeclareFieldFormat{eprint:customeprint}{%
	\em Available at \upshape \url{#1}}

\DeclareFieldFormat{eprint:online}{%
	Available \href{https://#1}{online}}

\DeclareSourcemap{
	\maps[datatype=bibtex]{
		\map{
			\step[fieldsource=mrnumber,		fieldtarget=eprint, final]
			\step[fieldset=eprinttype,		fieldvalue=mrnumber]
		}
		\map{
			\step[fieldsource=arxiv,		fieldtarget=eprint, final]
			\step[fieldset=eprinttype,		fieldvalue=arxiv]
		}
		\map{
			\step[fieldsource=jstor,		fieldtarget=eprint, final]
			\step[fieldset=eprinttype,		fieldvalue=jstor]
		}
		\map{
			\step[fieldsource=customeprint,	fieldtarget=eprint, final]
			\step[fieldset=eprinttype,		fieldvalue=customeprint]
		}
		\map{
			\step[fieldsource=online,		fieldtarget=eprint, final]
			\step[fieldset=eprinttype,		fieldvalue=online]
		}
	}
}

\DeclareFieldFormat{doi}{%
	\href{https://doi.org/#1}{DOI}%
}

\DeclareFieldFormat{url}{%
	\href{https://#1}{\texttt{#1}}%
}

\DeclareFieldFormat{comments}{%
	\emph{#1}%
}

\DeclareFieldFormat{manual}{%
	\emph{#1}%
} 
 
\DeclareFieldFormat
	[article,inbook,incollection,inproceedings,patent,thesis,unpublished]
	{title}{{#1\isdot}}

\DeclareFieldFormat
	[article,book,inbook,incollection,inproceedings,patent,thesis,unpublished, online]
	{date}{\mkbibparens{#1}}
\DeclareFieldFormat
	[article,book,inbook,incollection,inproceedings,patent,thesis,unpublished]
	{volume}{\mkbibbold{#1}}
\DeclareFieldFormat
	{pages}{\mkbibparens{#1}}

\DeclareFieldFormat{edition}{%
	\ifinteger{#1}
	{\mkbibordedition{#1}~\bibstring{edition}\addcomma}
	{#1~\bibstring{edition}\addcomma}}

\renewbibmacro*{volume+number+eid}{%
	\printfield{volume}%
\setunit*{\adddot}%
	\printfield{number}%
\setunit{\space}%
	\printfield{eid}%
}

\DeclareBibliographyDriver{article}{%
	\usebibmacro{bibindex}%
	\usebibmacro{begentry}%
	\usebibmacro{author}%
\setunit{\addspace}%
	\usebibmacro{date}%
\newunit\newblock
	\usebibmacro{title}%
\newunit\newblock
	\printfield{journaltitle}\addperiod%
\setunit*{\addspace}%
	\usebibmacro{volume+number+eid}%
\setunit*{\addspace}\newblock
	\printfield{pages}
\setunit*{\addnbspace}
	\printfield{comments}%
	\usebibmacro{eprint}%
\setunit*{\addnbspace}%
	\printfield{doi}%
	\usebibmacro{finentry}%
}

\DeclareBibliographyDriver{online}{%
	\usebibmacro{bibindex}%
	\usebibmacro{begentry}%
	\usebibmacro{author}%
\setunit{\addspace}%
	\usebibmacro{date}%
\newunit\newblock
	\usebibmacro{title}%
\newunit\newblock
	\printfield{comments}%
\setunit*{\addspace}
	\usebibmacro{eprint}%
	\usebibmacro{finentry}%
}

\DeclareBibliographyDriver{book}{%
	\usebibmacro{bibindex}%
	\usebibmacro{begentry}%
	\usebibmacro{author}%
\setunit{\addspace}\newblock
	\usebibmacro{date}%
\newunit\newblock
	\usebibmacro{title}%
\setunit*{\addspace}\newblock
	\printfield{edition}%
\newunit\newblock
	\printlist{publisher}
\setunit*{\addspace}%
	\printfield{volume}%
\setunit*{\addspace}\newblock%
	\printfield{comments}%
	\usebibmacro{eprint}%
\setunit*{\addspace}
	\printfield{doi}
	\usebibmacro{finentry}%
}

\DeclareBibliographyDriver{thesis}{%
	\usebibmacro{bibindex}%
	\usebibmacro{begentry}%
	\usebibmacro{author}%
\setunit{\addspace}\newblock
	\usebibmacro{date}%
\newunit\newblock
	\usebibmacro{title}.%
\setunit*{\addspace}\newblock
	\space Thesis, \printlist{institution}
	\usebibmacro{eprint}%
	\printfield{comments}%
\setunit*{\addspace}
	\printfield{doi}%
	\usebibmacro{finentry}%
}

\DeclareBibliographyDriver{inproceedings}{%
	\usebibmacro{bibindex}%
	\usebibmacro{begentry}%
	\usebibmacro{author}%
\setunit{\addspace}%
	\usebibmacro{date}%
\newunit\newblock
	\usebibmacro{title}%
\newunit\newblock
	\printfield{booktitle}\addcomma%
\setunit*{\addspace}%
	\printfield{series}\addcomma%
\setunit*{\addspace}\newblock
	\usebibmacro{editor}\addcomma
\setunit*{\addspace}\newblock
	\printlist{publisher}
\setunit*{\addspace}%
	\printfield{volume}%
\setunit*{\addspace}%
	\printfield{pages}
\setunit*{\addspace}
	\usebibmacro{eprint}%
	\printfield{comments}%
\setunit*{\addnnspace}
	\printfield{doi}
	\usebibmacro{finentry}%
}

\DeclareBibliographyDriver{inbook}{%
	\usebibmacro{bibindex}%
	\usebibmacro{begentry}%
	\usebibmacro{author}%
\setunit{\addspace}%
	\usebibmacro{date}%
\newunit\newblock
	\usebibmacro{title}%
\newunit\newblock
	\printfield{booktitle}\addcomma%
\setunit*{\addspace}%
	\printfield{series}\addcomma%
\setunit*{\addspace}\newblock
	\usebibmacro{editor}\addcomma
\setunit*{\addspace}\newblock
	\printlist{publisher}
\setunit*{\addspace}%
	\printfield{volume}%
\setunit*{\addspace}%
	\printfield{pages}
\setunit*{\addspace}
	\usebibmacro{eprint}%
	\printfield{comments}%
\setunit*{\addspace}
	\printfield{doi}
	\usebibmacro{finentry}%
}

\DeclareBibliographyDriver{incollection}{%
	\usebibmacro{bibindex}%
	\usebibmacro{begentry}%
	\usebibmacro{author}%
\setunit{\addspace}%
	\usebibmacro{date}%
\newunit\newblock
	\usebibmacro{title}%
\newunit\newblock
	\printfield{booktitle}\addcomma%
\setunit*{\addspace}%
	\printfield{series}\addcomma%
\setunit*{\addspace}\newblock
	\printlist{publisher}
\setunit*{\addspace}%
	\printfield{volume}%
\setunit*{\addspace}%
	\printfield{pages}
\setunit*{\addspace}
	\usebibmacro{eprint}%
	\printfield{comments}%
\setunit*{\addspace}
	\printfield{doi}
	\usebibmacro{finentry}%
}

\AtEveryBibitem{%
	\clearfield{day}%
	\clearfield{month}%
} 
 
\DeclareNameAlias{sortname}{given-family}
\DeclareNameAlias{default}{given-family}

\usepackage{geometry}
\usepackage{indentfirst, environ}
\usepackage[utf8]{inputenc}

\usepackage[UKenglish]{babel}
\usepackage{amsmath, amsthm, amssymb, mathrsfs, bm, mathtools}
\usepackage{wasysym}
\allowdisplaybreaks

\usepackage{titlesec}
\usepackage{xifthen}
\usepackage{xspace}
\usepackage{xparse}

\newcommand{\defac}[2]{(#1)_{#2}}

\newcommand{\RR}{R}

\newcommand{\ii}{i}
\newcommand{\II}{I}

\newcommand{\nn}{n}

\renewcommand{\tt}{t}

\newcommand{\cc}{c}

\renewcommand{\aa}{\alpha}
\newcommand{\jj}{j}
\newcommand{\mm}{m}
\newcommand{\MM}{M}
\newcommand{\pp}{p}

\newcommand{\rr}{r}
\newcommand{\xx}{x}
\newcommand{\yy}{y}
\newcommand{\zz}{z}

\newcommand{\triv}{\textnormal{triv}}

\usepackage{titlesec}
\titleformat{\section}
{\sffamily \Large \bfseries \bm}{\thesection}{1em}{}
\titleformat{\subsection}
{\sffamily \large \bfseries \bm}{\thesubsection}{1em}{}
\titleformat{\subsubsection}
{\sffamily \normalsize \bfseries \bm}{\thesubsubsection}{1em}{}
\titleformat{\paragraph}
{\sffamily \normalsize \bfseries \bm}{\theparagraph}{1em}{}
\titleformat{\subparagraph}[runin]
{\sffamily \normalsize \bfseries \bm}{\thesubparagraph}{1em}{}

\usepackage{enumitem}

\newcommand{\romannumbering}{%
	\renewcommand{\labelenumi}{\upshape(\roman{enumi})}%
	\renewcommand{\theenumi}{\upshape(\roman{enumi})}}

\usepackage{scrextend}

\setlist[description]{%
	topsep	= \smallskipamount,		
	itemsep	= \smallskipamount,		
	font	= {\mdseries\itshape},	
}

\newlist{ranlist}{enumerate}{3}
\setlist[ranlist,1]{
	label=(\roman*) }
\setlist[ranlist,2]{
	label=(\textit{\alph*}),
	ref=(\roman{ranlisti}.\textit{\alph*})	}
\setlist[ranlist,3]{
	label=\arabic*.,
	ref=(\roman{ranlisti}.\textit{\alph{ranlistii}}.\arabic*)	}

\usepackage[dvipsnames,hyperref]{xcolor}

\def\IfAmpersandUseAlign#1#2&#3\EndIfAmpersandUseAlign%
{%
	\if\relax\detokenize{#3}\relax
	\begin{equation*}%
	#1%
	\end{equation*}%
	\else
	\begin{align*}%
	#1%
	\end{align*}%
	\fi
}
\def\[#1\]%
{%
	\IfAmpersandUseAlign{#1}#1&\EndIfAmpersandUseAlign
}

\newcommand{\one}  [1]{\bm1( #1 )}

\let\originalexp\exp
\let\exp\relax
\DeclareRobustCommand{\exp} [1]{\originalexp(#1)}

\newcommand{\expb}[2][]{
	\ifthenelse{\equal{}{#1}}
	{\originalexp\bigl( #2 \bigr)}
	{\originalexp_{#1}\bigl( #2 \bigr)}
}

\newcommand{\logn}[1][]{
	\ifthenelse{\equal{}{#1}}
	{\log n}
	{(\log n)^{#1}}
}
\newcommand{\logk}[1][]{
	\ifthenelse{\equal{}{#1}}
	{\log k}
	{(\log k)^{#1}}
}

\newcommand{\Quad}[1]{
	\mathchoice
	{\quad\text{#1}\quad}
	{\text{ #1 }}
	{\text{ #1 }}
	{\text{ #1 }}
}

\newcommand{\Qand}{\Quad{and}}
\newcommand{\Qfor}{\Quad{for}}
\newcommand{\Qforall}{\Quad{for all}}
\newcommand{\Qwhere}{\Quad{where}}
\newcommand{\Qwith}{\Quad{with}}

\usepackage{calc}
\newlength{\halfplusheight}
\newcommand{\MAX}[1]{%
	\mathop{\raisebox{\halfplusheight}{\(\displaystyle\max_{#1}\)}}%
}

\newcommand{\LIM}[1]{%
	\mathop{\raisebox{\halfplusheight}{\(\displaystyle\lim_{#1}\)}}%
}
\newcommand{\LIMSUP}[1]{%
	\mathop{\raisebox{\halfplusheight}{\(\displaystyle\limsup_{#1}\)}}%
}
\newcommand{\LIMINF}[1]{%
	\mathop{\raisebox{\halfplusheight}{\(\displaystyle\liminf_{#1}\)}}%
}

\newcommand{\bcdot}{\ensuremath{\bm{\cdot}}}
\newcommand{\cq}{\coloneqq}

\newcommand{\abs}  [1]{\lvert #1 \rvert}
\newcommand{\absb} [1]{\big\lvert #1 \bigr\rvert}
\newcommand{\absB} [1]{\Bigl\lvert #1 \Bigr\rvert}
\newcommand{\absbb}[1]{\biggl\lvert #1 \biggr\rvert}

\newcommand{\rbr} [1]{ ( #1 ) }
\newcommand{\rbb} [1]{\bigl( #1 \bigr)}
\newcommand{\rbB} [1]{\Bigl( #1 \Bigr)}
\newcommand{\rbbb}[1]{\biggl( #1 \biggr)}

\newcommand{\sbb} [1]{\bigl[ #1 \bigr]}

\newcommand{\bra} [1]{ \{ #1 \} }
\newcommand{\brb} [1]{\bigl\{ #1 \bigr\}}

\newcommand{\midb}{\bigm|}

\newcommand{\Oh}  [1]{\mathcal{O}( #1 )}
\newcommand{\Ohb} [1]{\mathcal{O}\bigl( #1 \bigr)}

\newcommand{\oh}  [1]{o( #1 )}

\newcommand{\pr}[2][]{
	\mathchoice
	{\ifthenelse{\isempty{#1}}
		{\mathbb{P}\bigl(#2\bigr)}
		{\mathbb{P}_{#1}\bigl(#2\bigr)}}
	{\ifthenelse{\isempty{#1}}
		{\mathbb{P}(#2)}
		{\mathbb{P}_{#1}(#2)}}
	{\ifthenelse{\isempty{#1}}
		{\mathbb{P}(#2)}
		{\mathbb{P}_{#1}(#2)}}
	{\ifthenelse{\isempty{#1}}
		{\mathbb{P}(#2)}
		{\mathbb{P}_{#1}(#2)}}
}
\newcommand{\prb}[2][]{
	\ifthenelse{\equal{}{#1}}
	{\mathbb{P}\bigl( #2 \bigr)}
	{\mathbb{P}_{#1}\bigl( #2 \bigr)}
}
\newcommand{\prB}[2][]{
	\ifthenelse{\equal{}{#1}}
	{\mathbb{P}\Bigl( #2 \Bigr)}
	{\mathbb{P}_{#1}\Bigl( #2 \Bigr)}
}
\newcommand{\prbb}[2][]{
	\ifthenelse{\equal{}{#1}}
	{\mathbb{P}\biggl( #2 \biggr)}
	{\mathbb{P}_{#1}\biggl( #2 \biggr)}
}
\newcommand{\prBB}[2][]{
	\ifthenelse{\equal{}{#1}}
	{\mathbb{P}\Biggl( #2 \Biggr)}
	{\mathbb{P}_{#1}\Biggl( #2 \Biggr)}
}
\newcommand{\prs}[2][]{
	\ifthenelse{\equal{}{#1}}
	{\mathbb{P}\left( #2 \right)}
	{\mathbb{P}_{#1}\left( #2 \right)}
}

\newcommand{\ex}[2][]{
	\mathchoice
	{\ifthenelse{\isempty{#1}}
		{\mathbb{E}\bigl(#2\bigr)}
		{\mathbb{E}_{#1}\bigl(#2\bigr)}}
	{\ifthenelse{\isempty{#1}}
		{\mathbb{E}(#2)}
		{\mathbb{E}_{#1}(#2)}}
	{\ifthenelse{\isempty{#1}}
		{\mathbb{E}(#2)}
		{\mathbb{E}_{#1}(#2)}}
	{\ifthenelse{\isempty{#1}}
		{\mathbb{E}(#2)}
		{\mathbb{E}_{#1}(#2)}}
}
\newcommand{\exb}[2][]{
	\ifthenelse{\equal{}{#1}}
	{\mathbb{E}\bigl( #2 \bigr)}
	{\mathbb{E}_{#1}\bigr( #2 \bigr)}
}
\newcommand{\exB}[2][]{
	\ifthenelse{\equal{}{#1}}
	{\mathbb{E}\Bigl( #2 \Bigr)}
	{\mathbb{E}_{#1}\Bigl( #2 \Bigr)}
}
\newcommand{\exbb}[2][]{
	\ifthenelse{\equal{}{#1}}
	{\mathbb{E}\biggl( #2 \biggr)}
	{\mathbb{E}_{#1}\biggl( #2 \biggr)}
}
\newcommand{\exBB}[2][]{
	\ifthenelse{\equal{}{#1}}
	{\mathbb{E}\Biggl( #2 \Biggr)}
	{\mathbb{E}_{#1}\Biggl( #2 \Biggr)}
}

\newcommand{\Varb}[2][]{
	\ifthenelse{\equal{}{#1}}
	{\mathbb{V}\textnormal{ar} \bigl(#2\bigr)}
	{\mathbb{V}\textnormal{ar}_{#1} \bigl(#2\bigr)}
}
\newcommand{\VAR}[2][]{
	\ifthenelse{\equal{}{#1}}
	{\textnormal{Var}(#2)}
	{\textnormal{Var}_{#1}(#2)}
}

\newcommand{\maxt}[1]{
	\mathchoice
	{\textstyle \max_{#1} \displaystyle}
	{\max_{#1}}
	{\max_{#1}}
	{\max_{#1}}
}

\newcommand{\binomt}[2]{
	\mathchoice
	{\textstyle \binom{#1}{#2} \displaystyle}
	{\binom{#1}{#2}}
	{\binom{#1}{#2}}
	{\binom{#1}{#2}}
}

\newcommand{\sumt}[2][]{
	\mathchoice
	{\ifthenelse{\isempty{#1}}
		{\textstyle \sum_{#2}      \displaystyle}
		{\textstyle \sum_{#2}^{#1} \displaystyle}}
	{\ifthenelse{\isempty{#1}}
		{\sum_{#2}}
		{\sum_{#2}^{#1}}}
	{\ifthenelse{\isempty{#1}}
		{\sum_{#2}}
		{\sum_{#2}^{#1}}}
	{\ifthenelse{\isempty{#1}}
		{\sum_{#2}}
		{\sum_{#2}^{#1}}}
}
\newcommand{\sumd}[2][]{
	\ifthenelse{\isempty{#1}}
		{\sum_{#2}}
		{\sum_{#2}^{#1}}
}

\newcommand{\intt}[2][]{
	\mathchoice
	{\ifthenelse{\isempty{#1}}
		{\textstyle \int_{#2}      \displaystyle}
		{\textstyle \int_{#2}^{#1} \displaystyle}}
	{\ifthenelse{\isempty{#1}}
		{\int_{#2}}
		{\int_{#2}^{#1}}}
	{\ifthenelse{\isempty{#1}}
		{\int_{#2}}
		{\int_{#2}^{#1}}}
	{\ifthenelse{\isempty{#1}}
		{\int_{#2}}
		{\int_{#2}^{#1}}}
}

\newcommand{\prodt}[2][]{
	\mathchoice
	{\ifthenelse{\isempty{#1}}
		{\textstyle \prod_{#2}      \displaystyle}
		{\textstyle \prod_{#2}^{#1} \displaystyle}}
	{\ifthenelse{\isempty{#1}}
		{\prod_{#2}}
		{\prod_{#2}^{#1}}}
	{\ifthenelse{\isempty{#1}}
		{\prod_{#2}}
		{\prod_{#2}^{#1}}}
	{\ifthenelse{\isempty{#1}}
		{\prod_{#2}}
		{\prod_{#2}^{#1}}}	
}
\newcommand{\prodd}[2][]{
	\ifthenelse{\isempty{#1}}
		{\prod_{#2}}
		{\prod_{#2}^{#1}}
}

\newcommand{\toinf}[1]{\ensuremath{#1\to\infty}}
\newcommand{\asinf}[1]{\text{as \(#1\to\infty\)}}

\newcommand{\aszero}[1]{\text{as \(#1\to0\)}}

\DeclareMathOperator{\Unif}{Unif}

\DeclareMathOperator{\Pois}{Pois}

\DeclareMathOperator{\Bin}{Bin}

\DeclareMathOperator{\N}{N}

\DeclareMathOperator{\tr}{tr}

\DeclareMathOperator{\Fix}{Fix}
\newcommand{\symgr}{\mathcal S}
\newcommand{\altgr}{\mathcal A}

\newcommand{\TV}{\textnormal{TV}}

\newcommand{\typ}{\textnormal{typ}}
\newcommand{\id}{\textnormal{id}}

\newcommand{\tmix}{t_\mix}

\newcommand{\mix}{\textnormal{mix}}

\newcommand{\MT}{\textnormal{MT}}
\newcommand{\ET}{\textnormal{ET}}

\newcommand{\Ninn}{{N\in\mathbb{N}}}

\newcommand{\mbc}{\mathbb{C}}

\newcommand{\mbf}{\mathbb{F}}

\newcommand{\mbn}{\mathbb{N}}

\newcommand{\mbr}{\mathbb{R}}

\newcommand{\mbz}{\mathbb{Z}}

\newcommand{\mce}{\mathcal{E}}

\newcommand{\mcp}{\mathcal{P}}

\newcommand{\mcs}{\mathcal{S}}

\newcommand{\mcx}{\mathcal{X}}

\newcommand{\msf}{\mathscr{F}}
\newcommand{\msg}{\mathscr{G}}

\newcommand{\nt}{\addtocounter{equation}{1}\tag{\theequation}}

\newenvironment{Proof}[1][\proofname]{%
	\proof[\upshape\bfseries\sffamily\boldmath#1]
}{\endproof}

\newenvironment{center-small}
	{\par\centering\smallskip}
	{\par\smallskip}

\newcommand{\eps}{\varepsilon}

\usepackage{manyfoot}
\SetFootnoteHook{\hspace*{-1.8em}}
\DeclareNewFootnote{bl}[gobble]
\newcommand{\blfootnote}[1]{\footnotebl{\sffamily#1}}

\makeatletter
\newenvironment{subtheorem}[1]{%
	\def\subtheoremcounter{#1}%
	\refstepcounter{#1}%
	\protected@edef\theparentnumber{\csname the#1\endcsname}%
	\setcounter{parentnumber}{\value{#1}}%
	\setcounter{#1}{0}%
	\expandafter\def\csname the#1\endcsname{\theparentnumber\alph{#1}}%
	\expandafter\def\csname theH#1\endcsname{thm.\theparentnumber\alph{#1}}%
	\unskip\ignorespaces
}{%
	\setcounter{\subtheoremcounter}{\value{parentnumber}}%
	\ignorespacesafterend
}
\makeatother
\newcounter{parentnumber}

\makeatletter
\newenvironment{subtheorem-num}[1]{%
	\def\subtheoremcounter{#1}%
	\refstepcounter{#1}%
	\protected@edef\theparentnumber{\csname the#1\endcsname}%
	\setcounter{parentnumber}{\value{#1}}%
	\setcounter{#1}{0}%
	\expandafter\def\csname the#1\endcsname{\theparentnumber.\arabic{#1}}%
	\expandafter\def\csname theH#1\endcsname{thm.\theparentnumber.\arabic{#1}}%
	\unskip\ignorespaces
}{%
	\setcounter{\subtheoremcounter}{\value{parentnumber}}%
	\ignorespacesafterend
}
\makeatother

\usepackage[
	colorlinks,
	pdfauthor = {Evita\ Nestoridi\ and\ Sam\ Olesker-Taylor},
	pdfusetitle
]{hyperref}
\hypersetup{
	colorlinks,
	linkcolor	= {blue},
	urlcolor	= {blue},
	citecolor	= {ForestGreen}
}

\usepackage[capitalise]{cleveref}

\newtheoremstyle{sfsl}
{1\baselineskip}		
{1\baselineskip}		
{\slshape}				
{}						
{\bfseries\sffamily}	
{.}						
{0.5em}					
{\thmname{#1}\thmnumber{ #2}\thmnote{ \textnormal{\sffamily(#3)}}}

\newtheoremstyle{sfup}
{1\baselineskip}		
{1\baselineskip}		
{\upshape}				
{}						
{\bfseries\sffamily}	
{.}						
{0.5em}					
{\thmname{#1}\thmnumber{ #2}\thmnote{ \textnormal{\sffamily(#3)}}}

\theoremstyle{sfsl}

\newtheorem*{thm*}{Theorem}
\newtheorem{thm} {Theorem}[section]
\crefname{thm}{Theorem}{Theorems}

\newtheorem*{introthm*}{Theorem}
\newtheorem{introthm}{Theorem}

\crefname{introthm}{Theorem}{Theorems}

\newtheorem*{cor*}{Corollary}
\newtheorem{cor} [thm]{Corollary}
\crefname{cor}{Corollary}{Corollaries}

\newtheorem*{introcor*}{Corollary}
\newtheorem{introcor}{Corollary}

\crefname{introcor}{Corollary}{Corollaries}

\newtheorem*{lem*}    {Lemma}
\newtheorem{lem} [thm]{Lemma}
\crefname{lem}{Lemma}{Lemmas}

\newtheorem*{introlem*}{Lemma}
\newtheorem{introlem}{Lemma}

\crefname{introlem}{Lemma}{Lemmas}

\newtheorem*{prop*}    {Proposition}
\newtheorem{prop} [thm]{Proposition}
\crefname{prop}{Proposition}{Propositions}

\newtheorem*{clm*}    {Claim}

\crefname{clm}{Claim}{Claims}

\newtheorem*{defn*}    {Definition}
\newtheorem{defn} [thm]{Definition}
\crefname{defn}{Definition}{Definitions}

\newtheorem*{introdefn*}{Definition}
\newtheorem{introdefn}{Definition}

\crefname{introdefn}{Definition}{Definitions}

\newtheorem{innercustomgeneric}{\customgenericname}
\providecommand{\customgenericname}{}
\newcommand{\newcustomtheorem}[2]{%
	\newenvironment{#1}[1]
	{%
		\renewcommand\customgenericname{#2}%
		\renewcommand\theinnercustomgeneric{##1}%
		\innercustomgeneric
	}
	{\endinnercustomgeneric}
}

\newcustomtheorem{customthm}{Theorem}
\newcustomtheorem{customprop}{Proposition}
\newcustomtheorem{customlemma}{Lemma}
\newcustomtheorem{customconj}{Conjecture}
\newcustomtheorem{customquestion}{Question}

\newtheorem*{conj*}   {Conjecture}

\crefname{conj}{Conjecture}{Conjectures}

\newenvironment{conj-ind*}
	{\begin{quote}\textsf{\textbf{Conjecture.}}\slshape}
	{\end{quote}}
\newenvironment{conj-ind}
	{\begin{quote}\vspace{-\glueexpr\baselineskip+\topsep}\begin{customconj}}
	{\end{customconj}\end{quote}}

\newenvironment{question-ind*}
	{\begin{quote}\textsf{\textbf{Conjecture.}}\slshape}
	{\end{quote}}
\newenvironment{question-ind}
	{\begin{quote}\vspace{-\glueexpr\baselineskip+\topsep}\begin{customquestion}}
	{\end{customquestion}\end{quote}}

\newtheorem*{hypothesis*}{Hypothesis}

\newtheorem*{hyp*}{Hypotheses}

\crefname{hyp}{Hypotheses}{Hypotheses}

\theoremstyle{sfup}

\crefname{exm} {Example}{Examples}
\crefname{exmT}{Example}{Examples}

\newtheorem{rmkT}[thm]{Remark}
	\newenvironment{rmkt}
	{\pushQED{\qed}\rmkT}
	{\popQED\endrmkT}
\crefname{rmk} {Remark}{Remarks}
\crefname{rmkT}{Remark}{Remarks}

\newtheorem*{rmk*} {Remark}
\newtheorem*{rmkTT}{Remark}
\newenvironment{rmkt*}
	{\pushQED{\qed}\rmkTT}
	{\popQED\endrmkTT}

\crefname{rmks} {Remarks}{Remarks}
\crefname{rmksT}{Remarks}{Remarks}

\newtheorem*{rmks*} {Remarks}
\newtheorem*{rmksTT}{Remarks}
\newenvironment{rmkst*}
	{\pushQED{\qed}\rmksTT}
	{\popQED\endrmksTT}

\crefname{intrormk} {Remark}{Remarks}
\crefname{intrormkT}{Remark}{Remarks}

\newtheorem*{intrormk*} {Remark}
\newtheorem*{intrormkTT}{Remark}
\newenvironment{intrormkt*}
	{\pushQED{\qed}\intrormkTT}
	{\popQED\endintrormkTT}

\newtheorem*{exm*} {Example}
\newtheorem*{exmTT}{Lemma}
	\newenvironment{exmt*}
	{\pushQED{\qed}\exmTT}
	{\popQED\endexmTT}

\newtheorem*{note*} {Note}
\newtheorem*{noteTT}{Note}
	\newenvironment{notet*}
	{\pushQED{\qed}\noteTT}
	{\popQED\endnoteTT}

\numberwithin{equation}{section}

\title{\sffamily Limit Profiles for Reversible Markov Chains}
\author{\sffamily Evita Nestoridi\qquad Sam Olesker-Taylor}
\date{}

\begin{document}

\maketitle

\blfootnote{%
	\quad%
	Evita Nestoridi,\quad%
	\href{mailto:exn@princeton.edu}{exn@princeton.edu}
\hfill%
	\href{mailto:oleskertaylor.sam@gmail.com}{oleskertaylor.sam@gmail.com},\quad%
	Sam Olesker-Taylor%
	\quad{}
\\%
	Department of Mathematics, Princeton University, USA%
\hfill%
	Statistical Laboratory, University of Cambridge, UK%
\par\smallskip\par
\centering%
	This research was supported by the EPSRC:
	EN by EP/R022615/1
and
	SOT by Doctoral Training Grant 1885554%
}

\vspace{-6ex}

\renewcommand{\abstractname}{\sffamily Abstract}


\begin{abstract}\noindent
	In a recent breakthrough, Teyssier~\cite{T:limit-profile} introduced a new method for approximating the distance from equilibrium of a random walk on a group. He used it to study the limit profile for the random transpositions card shuffle. His techniques were restricted to conjugacy-invariant random walks on groups; we derive similar approximation lemmas for random walks on homogeneous spaces and for general reversible Markov chains. We illustrate applications of these lemmas to some famous problems: the $k$-cycle shuffle, sharpening results of Hough~\cite{H:cutoff-k-cycle} and Berestycki, Schramm and Zeitouni~\cite{BSZ:k-cycle}; the Ehrenfest urn diffusion with many urns, sharpening results of Ceccherini-Silberstein, Scarabotti and Tolli~\cite{CsST:gelfand-applications}; a Gibbs sampler, which is a fundamental tool in statistical physics, with Binomial prior and hypergeometric posterior, sharpening results of Diaconis, Khare and Saloff-Coste~\cite{DKSc:exp-families}.
\end{abstract}



\small
\begin{quote}
\begin{description}
	\item [Keywords:]
	cutoff, limit profiles, random walk on groups, symmetric group, representation theory, Fourier transform, characters, Gelfand pairs, homogeneous spaces, spherical functions, eigenvalues and eigenfunctions of Markov chains, spectral representations
	
	\item [MSC 2020 subject classifications:]
	20C15, 20C30; 43A30, 43A65, 43A90; 60B15, 60J10, 60J20
	
	\item [Acknowledgements:]
	We thank Persi Diaconis for helpful comments and suggestions
\end{description}
\end{quote}
\normalsize






%
\sffamily
\boldmath

\setcounter{tocdepth}{1}
\tableofcontents

\setcounter{tocdepth}{3}

\unboldmath
\normalfont

\romannumbering

\section{Introduction: TV Approximation Lemmas and Limit Profiles}

The cutoff phenomenon describes a situation where a Markov chain stays away from equilibrium for some time, but then converges to equilibrium very abruptly.
In rare cases, one can find an explicit function which describes this sharp transition, called the \textit{limit profile};
see, eg, \cite{BhS:cutoff-nonbacktracking,HOt:rcg:abe:cutoff,LP:ramanujan}.

In this paper, we develop a technique which allows us to well-approximate the distance from equilibrium, and hence study the limit profiles.
We consider the cases of
	general reversible Markov chains
	using a spectral decomposition
and
	random walks on homogeneous spaces, ie $X = G/K$ with $G$ a group and $K$ a subgroup of $G$
	using Fourier analysis.
The method is an extension of one introduced by \textcite{T:limit-profile} for random walks on Cayley graphs where the generating set is a union of conjugacy classes.
We then apply these techniques to prove the limit profile behaviours for
	the $k$-cycle shuffle,
	the multiple Ehrenfest urn model
and
	the Gibbs sampler with Binomial prior densities,
sharpening results of
\textcite{H:cutoff-k-cycle,BSZ:k-cycle,CsST:gelfand-applications,DKSc:exp-families}.


\subsection{Mixing Times and Limit Profiles}

Let $\Omega$ be a finite set and $P$ a transition matrix on $\Omega$.
Then $P^t(x,y)$ is the probability of moving from $x$ to $y$ in $t$ steps for all $x,y \in \Omega$ and all $t \in \mbn_0$.
If $P$ is irreducible and aperiodic, then the basic limit theorem of Markov chains tells us that $P^t(x,\cdot)$ converges to the (unique) invariant distribution $\pi$ as \toinf t with respect to the \textit{total variation} (abbreviated \textit{TV})~distance,
defined by%
\[
	d_\TV(\tt,\xx)
\cq
	d_\TV\rbr{ P^\tt(\xx,\cdot), \: \pi}
\cq
	\tfrac12 \sumt{\yy \in \Omega} \abs{P^\tt(\xx,\yy) - \pi(\yy)}
\Qfor
	\xx \in \Omega
\Qand
	\tt \in \mbn_0.
\]
The most common situation is to study the \textit{worst-case} TV distance:
\(
	d_\TV(\cdot)
\cq
	\maxt{\xx \in \Omega}
	d_\TV(\cdot, \xx).
\)
There are other possibilities, such as the \textit{typical} TV distance
where the starting point $\xx$ is chosen according to $\pi$:
\(
	d_\typ(\cdot)
\cq
	\sumt{\xx \in \Omega}
	\pi(\xx) d_\TV(\cdot, \xx).
\)
The (\textit{worst-case}) \textit{mixing time} is then defined by
\[
	\tmix(\eps)
\cq
	\inf\brb{\tt \ge 0 \midb d_\TV(\tt) \le \eps}
\Qfor
	\eps \in [0,1].
\]
For a sequence of Markov chains indexed by $N$, if there exist $(\tt_\star^{(N)})_\Ninn$ and $(w_\star^{(N)})_\Ninn$ satisfying
\[
	\LIM{\alpha \to -\infty}
	\LIMINF{\toinf N}
	d_\TV^{(N)}\rbb{ \tt_\star^{(N)} + \alpha w_\star^{(N)} }
=
	1
\Quad{and}
	\LIM{\alpha \to +\infty}
	\LIMSUP{\toinf N}
	d_\TV^{(N)}\rbb{ \tt_\star^{(N)} + \alpha w_\star^{(N)} }
=
	0,
\]
then the sequence of chains exhibits \textit{cutoff} at $t_\star$ with \textit{window}
$\Oh{w_*}$.

One can look beyond just finding the cutoff time and window, but instead determine the \textit{profile} inside the window:
	the aim is to choose $\tt_\star$ and $w_\star$ appropriately so that
	\[
		\varphi(\alpha)
	\cq
		\LIM{N \to \infty} d_\TV^{(N)}\rbb{ \tt_\star^{(N)} + \alpha w_\star^{(N)} }
	\quad
		\text{exists for all $\alpha \in \mbr$}.
	\]
The limit \toinf N is taken for each fixed $\alpha \in \mbr$.

Officially, when we look at $d_\TV(t)$, we need $t \in \mbn$; in practice, we omit floor/ceiling signs.

\subsection{TV Convergence Profile for Random Walks}

In this paper we present three lemmas for obtaining the TV profile for random walks; see \cref{res:tv-approx:rev,res:tv-approx:tey,res:tv-approx:hom}.
They work by finding a decomposition of the TV distance as a sum using either a spectral decomposition or Fourier analysis.
One then separates out the `important' terms in the sum to give a `main term' (which asymptotically captures all the TV mass) and an `error' term.
\cref{res:tv-approx:rev,res:tv-approx:hom} are original contributions; \cref{res:tv-approx:tey} is due to \textcite{T:limit-profile}.
For each lemma, we give an example application, establishing a limit profile of the TV convergence to equilibrium.

We denote the cdf of the standard normal distribution by $\Phi$ throughout the paper.

\subsubsection{Reversible Markov Chains}
\label{sec:intro:res:rev}

First we consider general reversible Markov chains on an arbitrary set $\Omega$.
The following lemma is based off the well-known spectral decomposition for a reversible Markov chain $P$:
\[
	P^\tt(\xx,\yy)
=
	\sumt[\abs \Omega]{\ii=1}
	\sumt{\yy \in \Omega}
	\pi(\yy) f_\ii(\xx) f_\ii(\yy) \lambda_\ii^\tt
\Qforall
	\xx,\yy \in \Omega
\Qand
	\tt \in \mbn_0,
\]
where
	$P^\tt(\xx,\yy)$ is the probability of moving from $\xx$ to $\yy$ in $\tt$ steps,
	$\pi$ is the invariant distribution
and
	$\bra{f_\ii, \lambda_\ii}_{\ii=1}^{\abs \Omega}$ are the eigenstatistics;
see \cite[Lemma~12.2]{LPW:markov-mixing}.
Recall that, for $\xx \in \Omega$ and $\tt \in \mbn_0$, we write $d_\TV(\tt,\xx)$ for the TV distance from $\pi$ after $t$ steps when started from $\xx$.

We come to our first contribution: the TV-approximation lemma for reversible Markov chains.

\begin{introlem}[Reversible Markov Chains]
\label{res:tv-approx:rev}
	Consider a reversible, irreducible and aperiodic Markov chain on a finite set $\Omega$ with invariant distribution $\pi$.
	Denote
		by
		\(
			-1 < \lambda_{\abs \Omega} \le \ldots \le \lambda_2 < \lambda_1 = 1
		\)
		its eigenvalues
	and
		by $f_{\abs \Omega}, ..., f_1$ its corresponding orthonormal (with respect to $\pi$) eigenvectors.
	For $\tt \in \mbn_0$ and $x \in \Omega$,
	denote by $d_\TV(t,x)$ the TV distance from equilibrium (ie $\pi$) of the Markov chain started from $x$.
	
	For
		all $\tt \in \mbn_0$,
		all $x \in \Omega$
	and
		all $\II \subseteq \bra{2, ..., \abs \Omega}$,
	we have
	\[
	\absb{
		d_\TV\rbr{\tt,x}
	-	\tfrac12 \sumt{y \in \Omega} \pi(y) \absb{ \sumt{\ii \in \II} f_\ii(x) f_\ii(y) \lambda_\ii^\tt }
	}
	\le
		\tfrac12 \sumt{\ii \notin \II} \abs{ f_\ii(x) } \abs{ \lambda_\ii }^\tt.
	\]
\end{introlem}

As an application of \cref{res:tv-approx:rev}, we determine the limit profile for a specific two-component Gibbs sampler, which is an important tool in statistical physics as explained in \cite[\S1]{DKSc:exp-families}.

Let $(\mcx, \msf, \mu)$ and $(\Theta, \msg, \pi)$ be two probability spaces.
The probability measure $\pi$ is called the \textit{prior}.
Let $\bra{ f_\theta(\cdot) }_{\theta \in \Theta}$ be a family of probability densities on $\mcx$ with respect to $\mu$.
These define a probability measure $\Pr$ on $\mcx \times \Theta$ by
\[
	\Pr(A \times B)
\cq
	\intt{B} \intt{A} f_\theta(x) d\mu(x) \, d\pi(\theta)
\Qfor
	(A, B) \in \msf \times \msg.
\]
The marginal density on $\mcx$ is given by
\(
	m(x)
\cq
	\intt{\theta} f_\theta(x) d\pi(\theta)
\)
for $x \in \mcx$.
The \textit{posterior} density with respect to the prior $\pi$ is defined by
\(
	\pi(\theta \vert x)
\cq
	f_\theta(x) / m(x)
\)
for $(x, \theta) \in \mcx \times \Theta$.

The ($\mcx$-chain) Gibbs sampler is defined informally as follows (each draw is independent):
\begin{center-small}
	\bcdot\ input $x$;\quad%
	\bcdot\ draw $\theta \sim \pi(\cdot \vert x)$;\quad%
	\bcdot\ draw $x' \sim f_\theta(\cdot)$;\quad%
	\bcdot\ output $x'$.%
\end{center-small}\noindent
Formally, it is the Markov chain defined by the transition kernel $P$ given by
\[
	P(x,x')
\cq
	\intt{\Theta}
	\pi(\theta \vert x) f_\theta(x') \, d\pi(\theta)
=
	\intt{\Theta}
	f_\theta(x) f_\theta(x') / m(x) \, d\pi(\theta)
\Qfor
	x,x' \in \mcx.
\]
Observe that $P$ is reversible with respect to $m$, ie the marginal density on $\mcx$.

\smallskip

We consider the special case of \textit{location families}:
\(
	f_\theta(x)
=
	g(x - \theta)
\)
for all
\(
	(x, \theta) \in \mcx \times \Theta.
\)
for some function $g$;
see \cite[\S5]{DKSc:exp-families}.
The Gibbs sampler can then be realised in the following way:
\begin{center-small}
	\bcdot\ input $x$;\quad%
	\bcdot\ draw $\theta \sim \pi(\cdot \vert x)$;\quad%
	\bcdot\ draw $\eps \sim g$;\quad%
	\bcdot\ output $x' \cq \theta + \eps$.%
\end{center-small}\noindent
We consider prior $\pi$ and $g$ each being Binomial, which leads to a hypergeometric posterior.

\smallskip

Our next contribution is
	the limit profile for the two-component Gibbs sampler with Binomial priors,
	established as an application of \cref{res:tv-approx:rev}.
A more refined statement is given in \cref{res:rev:gibbs:res}.

\begin{introthm}[Gibbs Sampler]
\label{res:intro:gibbs}
	Let $\nn_1, \nn_2 \in \mbn$ and $p \in (0,1)$;
	write $\nn \cq \nn_1 + \nn_2$ and $\aa \cq p/(1-p)$.
	Let $\pi \sim \Bin(\nn_1, p)$ and $g \sim \Bin(\nn_2, p)$.
	For $\tt \in \mbn_0$, write $d_\TV^{\nn_1,\nn_2,p}(t)$ for the TV distance of the (location family) Gibbs sampler after $\tt$ steps started from $0 \in \mbn$ from its invariant distribution $m$.
	
	Suppose that $\min\bra{p, 1-p} \cdot \nn \gg 1$.
	Then,
	for all $\cc \in \mbr$ (independent of $\nn$),
	we have
	\[
		d_\TV^{\nn_1,\nn_2,p}\rbb{
			\rbb{ \tfrac12 \log(\aa \nn) + \cc } / \log\rbb{ \tfrac1{1 - \nn_2/\nn} }
		}
	\to
		2 \, \Phi\rbb{ \tfrac12 e^{-c} } - 1.
	\]
\end{introthm}

The above set-up implicitly sets
	the sample spaces
		$\mcx \cq \bra{0, ..., \nn_1}$
	and
		$\Theta \cq \bra{0, ..., \nn_2}$
and
	the event spaces
		to be the respective set of all subsets.
The sample spaces are finite, so this is natural.

Cutoff for the $L_2$ mixing time of this Gibbs sampler was established by \textcite[\S5.1]{DKSc:exp-families}; these tools could likely be adapted to give cutoff for the usual TV ($L_1$) mixing time. However, the techniques of \textcite[\S5.1]{DKSc:exp-families} are not sufficiently refined to give access to the limit profile; a more detailed analysis is required.

\subsubsection{Random Walks on Groups}

We start by recalling some standard terminology from representation theory.

\begin{introdefn*}
	Let $G$ be a finite group and $V$ a finite dimensional vector space over $\mbc$.
	A \textit{representation} $\rho$ of $G$ over $V$ is an action
	\(
		(g,v) \mapsto \rho(g) \cdot v : G \times V \to V
	\)
	such that $\rho(g) : V \to V$ is an invertible linear map for all $g \in G$.
%
	The \textit{Fourier transform} of a function $\mu : G \mapsto \mbc$ with respect to the representation $(\rho, V)$ is the linear operator
	\(
		\widehat \mu(\rho) : V \to V
	\)
	defined by
	\(
		\widehat \mu(\rho)
	\cq
		\sumt{g \in G}
		\mu(g) \rho(g).
	\)
\end{introdefn*}

Using the Fourier inversion formula,
for all probability measures $\mu$ on $G$ and all $t \in \mbn_0$,
we~have
\[
	d_\TV(\mu^{*t}, \: \Unif_G)
=
	\tfrac12 \abs G^{-1}
	\sumt{g \in G}
	\absb{
		\sumt{\rho \in R^*}
		d_\rho \tr\rbb{ \widehat \mu(\rho)^t \rho(g^{-1}) }
	},
\]
where
	$\mu^{*t}$ is the $t$-fold self-convolution of $\mu$,
	$R^*$ is the set of all non-constant irreducible representations (abbreviated \textit{irreps}) of $G$
and
	$d_\rho$ is the dimension of the irrep~$\rho$;
see \cite[\S3.10]{CsST:harmonic-analysis-finite-groups}.

If $\mu$ is the step distribution of a random walk on $G$, then this determines exactly TV distance after $t$ steps; cf the well-known spectral representation for reversible random walks.
One must still control the Fourier transform at arbitrary irreps.
There are two important special cases.
\begin{itemize}[itemsep = 0pt, topsep = \smallskipamount, label = \bcdot]
	\item 
	Suppose that $\mu$ is \textit{conjugacy-invariant}, ie $\mu(g) = \mu(h^{-1} g h)$ for all $g,h \in G$.
	By Schur's lemma, $\widehat \mu(\rho)$ is a multiple of the identity for each irrep $\rho$.
	Then the key object in calculating the Fourier transform is the \textit{character}:
	\(
		\chi_g(\rho) \cq \tr(\rho(g))
	\)
	for $g \in G$ and $\rho \in R^*$.
	This is the case considered originally
	in \cite{DS:random-trans}, and then in \cite{T:limit-profile},
	for random transpositions.
	
	\item 
	Suppose that the matrices $\widehat \mu(\rho)$ have only one non-trivial entry which is in the first position
	(in an appropriate `spherical' basis).
	This radical but frequent simplification occurs in the framework of \textit{Gelfand pairs}; see \S\ref{sec:hom} for details.
	\textcite{DS:bernoulli-laplace} consider this in the set-up of the Bernoulli--Laplace urn model, and more generally.
\end{itemize}

\paragraph{Conjugacy-Invariant Random Walks}

In this subsection we state \citeauthor{T:limit-profile}'s lemma for conjugacy-invariant random walks.

\setcounter{introdefn}{\value{introlem}}
\begin{introdefn}
	A random walk on $G$ is \textit{conjugacy-invariant} if there is a probability measure $\mu$ which is constant on each conjugacy class of $G$ for which the transition matrix $P$ satisfies $P(x,xg) = \mu(g)$ for all $x \in G$.
	For a representation $\rho$, define the \textit{character ratio}
	\(
		s_\rho \cq d_\rho^{-1} \sumt{g \in G} \mu(g) \chi_g(\rho).
	\)
\end{introdefn}


\citeauthor{T:limit-profile}'s lemma for conjugacy-invariant random walks states the following.

\begin{introlem}[{\textcite[Lemma~2.1]{T:limit-profile}}]
\label{res:tv-approx:tey}
	Let $G$ be a finite group;
	let $\mu : G \mapsto [0,1]$ be a conjugacy-invariant probability distribution on $G$.
	For $\tt \in \mbn_0$, denote by $d_\TV(t)$ the TV distance to equilibrium of the random walk on $G$ started from the identity with step distribution $\mu$ and run for $\tt$ steps.

	Let $\tt \in \mbn_0$ and $I \subseteq \RR^*$, ie the set of non-trivial irreps of $G$.
	Then
	\(
		d_\TV(t) = d_\TV(\mu^{*t}, \: \Unif_G)
	\)~%
	and
	\[
	\absb{
		d_\TV(\tt)
		-	\tfrac12 \abs G^{-1} \sumt{g \in G} \absb{ \sumt{\rho \in I} d_\rho s_\rho^\tt \chi_\rho(g) }
	}
	\le
		\tfrac12 \sumt{\rho \in \RR^* \setminus I} d_\rho \abs{s_\rho}^\tt.
	\]
\end{introlem}

We apply this lemma to the $k$-cycle random walk on the symmetric group $\symgr_n$.
In this walk, at each step a $k$-cycle is chosen uniformly at random and composed with the current location.
We establish the limit profile for $2 \le k \ll n$.
There are parity constraints.
To handle such parity constraints, we follow the set-up used by \citeauthor{H:cutoff-k-cycle}:
\begin{itemize}[noitemsep, topsep = \smallskipamount, label = \bcdot]
	\item 
	if $k$ is odd, then the walk is supported on the set of even permutations;
	
	\item 
	if $k$ is even and $t$ is even, then the walk at time $t$ is supported on the set of even permutations;
	
	\item 
	if $k$ is even and $t$ is odd, then the walk at time $t$ is supported on the set of odd permutations.
\end{itemize}

We come to our next contribution:
	the limit profile for the random $k$-cycle shuffle,
	established as an application of \cref{res:tv-approx:tey}.
A more refined statement is given in \cref{res:k-cycle:res}.

\begin{introthm}[Random $k$-Cycles]
\label{res:intro:k-cycle}
	Let $k, n \in \mbn$ with $2 \le k \le n$.
	For $t \in \mbn_0$, denote by $d_\TV^{n,k}(t)$ the TV distance of the $k$-cycle random walk on $\symgr_n$ from the uniform distribution on the appropriate set of permutations of a fixed parity started from the identity and run for $t$ steps.
	
	Suppose that $2 \le k \ll n$.
	Then,
	for all $\cc \in \mbr$ (independent of $n$),
	we have
	\[
		d_\TV^{n,k}\rbb{ - n \rbr{ \log n + \cc } / \log(1 - k/n) }
	\to
		d_\TV\rbb{ \Pois(1 + e^{-\cc}), \: \Pois(1) }.
	\]
	If, further, $k \ll n / \log n$, then
	the same \TV\ limit holds when evaluated at time $\tfrac nk (\log n + \cc)$ instead.%
\end{introthm}

Cutoff for this shuffle was already been established by \textcite{H:cutoff-k-cycle}, for any $2 \le k \ll n$, using representation theory.
He also found the correct order of the window when $2 \le k \ll n / \log n$.
We handle any $2 \le k \ll n$, and find the \emph{precise} limit profile, not just the order of the window.

The case of random transpositions, ie $k = 2$, was one of the first Markov chains studied using representation theory; cutoff was established by \textcite{DS:random-trans}.
\textcite{BSZ:k-cycle} established cutoff for $k$ fixed, independent of $n$, using probabilistic arguments instead of representation theory.
\textcite{BS:cutoff-conj-inv} studied a generalisation where one draws uniformly from a prescribed conjugacy class with support $k$ with $2 \le k \ll n$.

The limit profile, even for $k = 2$, remained a famous open problem for a long time.
A breakthrough came recently by \textcite{T:limit-profile}, using \cref{res:tv-approx:tey} above;
we apply this lemma here.
Also, we adapt and extend some character theory for the $k$-cycle walk developed by \textcite{H:cutoff-k-cycle}.
Finally, we adapt and extend some of the analysis of \textcite{T:limit-profile} from $k = 2$ to general $k$.

\paragraph{Random Walks on Homogenous Spaces}

Finally we turn our attention to random walks on \textit{homogeneous spaces} $X \cq G/K$, where $G$ is a finite group and $K$ a subgroup of $G$.
Where \cite{DS:random-trans,T:limit-profile} considered conjugacy-invariant $\mu$ to simplify the calculation of the Fourier transforms,
here we consider the case
	that $\mu$ is $K$ bi-invariant, ie
	\(
		\mu(k_1 g k_2) = \mu(g)
	\)
	for all $g \in G$ and all $k_1, k_2 \in K$
and
	that $(G,K)$ is a \textit{Gelfand pair}, ie the algebra of $K$ bi-invariant functions (under convolution) is commutative;
	see \cref{def:hom:gelfand}.
In this case,
for any $K$ bi-invariant function $\mu$ on $G$,
if $(\rho, V)$ is a \textit{spherical} irrep, defined in \cref{def:hom:spherical:fn-rep}, then
the matrix $\widehat \mu(\rho)$ has only one non-zero entry, which is in the top-left position; this entry is called the \textit{spherical Fourier transform} of $\mu$ with respect to $\rho$ (rescaled by $\abs K$).
Moreover, if $(\tau, W)$ is a non-spherical irrep, then $\widehat \mu(\tau) = 0$ is the zero matrix.

Using this simplification, we prove the following lemma for random walks on homogeneous spaces corresponding to a Gelfand pair started from some element $\bar x \in K$ stabilised by $K$, ie $k \bar x = \bar x$ for all $k \in K$ (under the usual left coset action).
The canonical quotient projection $G \to G/K$ preserves the uniform distribution.
So the invariant distribution of any random walk on a homogenous space is uniform on that space.

Our next contribution is a TV-approximation lemma for random walks on homogeneous spaces.

\begin{introlem}[Homogeneous Spaces]
\label{res:tv-approx:hom}
	Let $(G, K)$ be a Gelfand pair and denote $X \cq G/K$.
	Let $\bar x$ be an element of $X$ whose stabiliser is $K$.
	Let
		$\bra{\varphi_i}_{i=0}^N$ be the associated spherical functions,
		with $\varphi_0(x) = 1$ for all $x \in X$,
		considered as $K$-invariant functions on $X$,
	and
		$\bra{d_i}_{i=0}^N$ the associated dimensions.
	Let $P$ be a $G$-invariant stochastic matrix and set $\mu_{\bar x}(\cdot) \cq P\rbr{\bar x, \cdot}$.
	For $\tt \in \mbn_0$, denote by $d_\TV(\tt,\bar x)$ the TV distance to equilibrium of the random walk on $X$ started from $\bar x$ with step distribution $\mu_{\bar x}$ and run for $\tt$ steps.
	
	Let $\tt \in \mbn_0$ and $I \subseteq \bra{1, ..., N}$.
	Then
	\(
		d_\TV(\tt, \bar x)
	=
		d_\TV\rbr{ \mu_{\bar x}^{*\tt}, \: \Unif_X }
	\)
	and
	\[
		\absB{
			d_\TV(\tt, \bar x)
		-	\tfrac12 \abs X^{-1} \sumt{x \in X} \absb{ \sumt{i \in I} d_i \varphi_i(x) \widehat \mu_{\bar x}(i)^\tt }
		}
	\le
		\tfrac12 \sumt{i \notin I} \sqrt{d_i} \abs{ \widehat \mu_{\bar x}(i) }^\tt,
	\]
	where $\widehat \mu_{\bar x} : i \mapsto \sumt{x \in X} \mu_{\bar x}(x) \overline{\varphi_i(x)}$ is the spherical Fourier transform of $\mu$ with respect to $\bra{\varphi_i}_{i=0}^N$.
\end{introlem}

We come to our final contribution:
	the limit profile for the multiple urn Ehrenfest urn diffusion model,
	established as an application of \cref{res:tv-approx:hom}.
A more refined statement is given in \cref{res:hom:e-urn:res}.

\begin{introthm}[Ehrenfest Urn]
\label{res:intro:e-urn}
	Let $n,m \in \mbn$.
	Consider $n$ labelled balls and $m+1$ labelled urns.
	Consider the following Markov chain:
		at each step,
		choose a ball and an urn uniformly and independently;
		place said ball in said urn.
	For $\tt \in \mbn_0$,
	denote by $d_\TV^{n,m}(\tt)$ the TV distance of this urn model started with all balls in a single urn from its invariant distribution and~run~for~$\tt$~steps.
	
	Suppose that $1 \le m \ll n$.
	Then,
		all $\cc \in \mbr$ (independent of $n$),
	we have
	\[
		d_\TV^{n,m}\rbb{ \tfrac12 n \log(nm) + \cc n }
	\to
		2 \, \Phi\rbb{ \tfrac12 e^{-\cc} } - 1.
	\]
\end{introthm}

Cutoff, but not the limit profile, was established for this multiple urn model by \textcite[\S6]{CsST:gelfand-applications} using representation theory.
To establish the profile,
we apply the approximation lemma for random walks on homogeneous spaces, ie Lemma \ref{res:tv-approx:hom}, using the character theory developed by \textcite{CsST:gelfand-applications}.

This model was originally introduced (with two urns) by \textcite{EE:e-urn} in \citeyear{EE:e-urn}.
In this case, the model can be viewed as a TV-preserving projection of the simple random walk on the $n$-hypercube.
There cutoff was established by \textcite[Example~3.19]{A:rw-group-mixing}.
The limit profile is even known:
see
	\textcite[Theorem~18 in \S6.2]{S:mixing-notes} (in French) for a `probabilistic' argument using convergence theorems
or
	\textcite[Theorem~1]{DGM:rw-hypercube} for a Fourier analytical argument.
We present a significantly simpler Fourier analytical argument, using only basic representation theory of the Abelian group $\mbz_2^d$
in \cref{res:rw-hypercube:res}.
\newpage

\subsubsection{Corollaries to TV Approximation Lemma for Reversible Markov Chains (\cref{res:tv-approx:rev})}

We close this section with two simple corollaries of the general TV-approximation lemma for reversible Markov chains, \cref{res:tv-approx:rev}.
The first is for transitive Markov chains; the second is for typical TV distance.
For transitive chains, the starting point is irrelevant;
that is, for each $t$, the map $x \mapsto d_\TV(t,x)$ is constant (ie does not depend on the input $x$).
In particular,
\(
	d_\TV(\cdot) = \sumt{x \in \Omega} \pi(x) d_\TV(\cdot,x).
\)
Also, by transitivity, the invariant distribution $\pi$ is uniform on $\Omega$.

\begin{subtheorem-num}{introcor}
	\label{res:tv-approx:rev:cor}

\begin{introcor}
\label{res:tv-approx:rev:cor:trans}
	Consider the set-up of \cref{res:tv-approx:rev};
	in addition, assume that the chain is transitive.
	
	For
		all $\tt \in \mbn_0$
	and
		$\II \subseteq \bra{2, ..., \abs \Omega}$,
	we have
	\[
		\absb{ d_\TV(\tt) - \tfrac12 \abs \Omega^{-2} \sumt{x,y \in \Omega} \absb{ \sumt{\ii \in \II} f_\ii(x) f_\ii(y) \lambda_\ii^\tt } }
	\le
		\tfrac12 \sumt{\ii \notin \II} \abs{ \lambda_\ii^\tt }.
	\]
\end{introcor}

Instead of looking at TV from a given starting point, we can also consider averaging over the starting point (with respect to the invariant distribution).
This is sometimes known as \textit{typical} TV distance (as opposed to \textit{worst-case}).
For $\tt \in \mbn_0$, denote
\[
	d_\typ(\tt) \cq \sumt{x \in \Omega} \pi(x) d_\TV(\cdot,x).
\]

\begin{introcor}
\label{res:tv-approx:rev:cor:typ}
	Consider the set-up of \cref{res:tv-approx:rev};
	no transitivity is necessary.
	
	For
		all $\tt \in \mbn_0$
	and
	 	all $\II \subseteq \bra{2, ..., \abs{\Omega}}$,
 	we have
	\[
		\absb{ d_\typ(\tt) - \tfrac12 \sumt{x,y \in \Omega} \pi(x) \pi(y) \absb{ \sumt{\ii \in \II} f_\ii(x) f_\ii(y) \lambda_\ii^\tt } }
	\le
		\tfrac12 \sumt{\ii \notin \II} \abs{ \lambda_\ii^\tt }.
	\]
\end{introcor}

\end{subtheorem-num}

\subsection{Organisation of the Paper}
\label{sec:intro:organisation}

The remainder of the paper is organised as follows.

\begin{itemize}[itemsep = 0pt, topsep = \smallskipamount, label = \bcdot]
	\item [\S\ref{sec:rev}]
	Here we study general reversible Markov chains.
	We prove the our TV-approximation lemma (\cref{res:tv-approx:rev}) via an application of the spectral decomposition for reversible Markov chains.
	
	As an application of \cref{res:tv-approx:rev},
	we establish the limit profile for a two-component Gibbs sampler, which are fundamental tools in statistical physics (see \cite[\S1]{DKSc:exp-families}).

	\item [\S\ref{sec:k-cycle}]
	Here we establish the limit profile of the random $k$-cycle walk on the symmetric group.
	We do this via an application of the TV-approximation lemma of \textcite{T:limit-profile} (\cref{res:tv-approx:tey}), along with extending and applying character theory developed by \textcite{H:cutoff-k-cycle}.
	
	\item [\S\ref{sec:hom}]
	Here we study random walks on homogeneous spaces corresponding to Gelfand pairs.
	We develop and apply (mostly classical) theory to prove our TV-approximation lemma (\cref{res:tv-approx:hom}).
	
	As an application of \cref{res:tv-approx:hom},
	we establish the limit profile for the famous Ehrenfest urn diffusion with many urns, using some character theory developed by \textcite{CsST:gelfand-applications}.
\end{itemize}


%

\section{Reversible Markov Chains}
\label{sec:rev}

In this section general reversible Markov chains are considered.
First we prove the lemma and corollaries from the introduction, then we apply them to a Gibbs sampler.

\subsection{Proof of TV-Approximation Lemmas for Reversible Markov Chains}

\cref{res:tv-approx:rev} follows from the usual spectral representation of TV distance along with some algebraic manipulations and inequalities.
\cref{res:tv-approx:rev:cor:typ,res:tv-approx:rev:cor:trans} follow, in an identical way to each other, from averaging both sides of \cref{res:tv-approx:rev} with respect to~$\pi$.
We give the full details now.

\begin{Proof}[Proof of \cref{res:tv-approx:rev}]
As an immediate consequence of \cite[Lemma~12.2]{LPW:markov-mixing},
for $x \in \Omega$,
we have
\[
	d_\TV(\tt,x)
=
	\tfrac12 \sumt{y \in \Omega} \pi(y)
	\absb{ \sumt[\abs \Omega]{\ii=2} f_\ii(x) f_\ii(y) \lambda_\ii^\tt }.
\]
Let $\II \subseteq \bra{2, ..., \abs \Omega }$.
Elementary manipulations using
	the triangle inequality (twice)
	then
		Cauchy--Schwarz
	and
		the fact that the eigenfunctions are orthonormal with respect to $\pi$,
give
\[
	\absb{
		d_\TV(\tt,x)
	-	\tfrac12 \sumt{y \in \Omega} \pi(y) \abs{ \sumt{\ii \in \II} f_\ii(x) f_\ii(y) \lambda_\ii^\tt }
	}
&
\le
	\tfrac12 \sumt{y \in \Omega} \pi(y) \abs{ \sumt{\ii \notin \II} f_\ii(x) f_\ii(y) \lambda_\ii^\tt }
\\&
\hspace{-7em}
\le
	\tfrac12 \sumt{\ii \notin \II}
	\abs{ f_\ii(x) \lambda_\ii^\tt }
	\sumt{y \in \Omega} \pi(y) \abs{ f_\ii(y) }
\le
	\tfrac12 \sumt{\ii \notin \II} \abs{ f_\ii(x) \lambda_\ii^\tt }.
\qedhere
\]
\end{Proof}

\begin{Proof}[Proof of \cref{res:tv-approx:rev:cor:typ,res:tv-approx:rev:cor:trans}]
For a transitive chain, for each $\tt$, the map $x \mapsto d_\TV(t,x)$ is constant (ie does not depend on the input $x$).
So we may replace $d_\TV(t,x)$ by $\sumt{x \in \Omega} \pi(x) d_\TV(t,x) = d_\typ(t)$.
The corollaries now follow by averaging the error term with respect to $\pi$, using Cauchy--Schwarz and the normalisation of the eigenfunctions.
\end{Proof}

%

\subsection{Application to Gibbs Sampler with Binomial Priors}
\label{sec:rev:gibbs}

\newcommand{\Kp}[2] {K_{#1}\rbr{#2; \pp, \nn}}
\newcommand{\Kpb}[2]{K_{#1}\rbb{#2; \pp, \nn}}
\newcommand{\Ka}[2] {K_{#1}\rbr{#2; \tfrac{\aa}{\aa+1}, \nn}}
\newcommand{\Kab}[2]{K_{#1}\rbb{#2; \tfrac{\aa}{\aa+1}, \nn}}

In this subsection,
we consider the Gibbs sampler with Binomial priors, namely $\pi \sim \Bin(\nn_1, \pp)$ and $g \sim \Bin(\nn_2, \pp)$, as described in \cref{res:intro:gibbs}.
Here $\mcx \cq \bra{0, 1, ..., \nn}$ where $\nn \cq \nn_1 + \nn_2$.

The following theorem is a restatement of \cref{res:intro:gibbs}, but written more formally: cutoff is for a sequence of Markov chains; we make this sequence explicit.

\begin{thm}
\label{res:rev:gibbs:res}
	Let $\nn_1, \nn_2 \in \mbn$ and $\pp \in (0,1)$;
	write $\nn \cq \nn_1 + \nn_2$ and $\aa \cq \pp / (1 - \pp)$.
	Consider the (location family) Gibbs sampler with $\pi \sim \Bin(\nn_1, \pp)$ and $g \sim \Bin(\nn_2, \pp)$.
	For $\tt \in \mbn_0$,
	let $d_\TV^{\nn_1, \, \nn_2, \, \pp}(\tt)$ denote the TV distance from equilibrium after $\tt$ steps in this Gibbs sampler started from 0.
	
	Let $(\nn_{1,N})_\Ninn, (\nn_{2,N})_\Ninn \in \mbn^\mbn$ and $(\pp_N) \in (0,1)^\mbn$;
	for each $\Ninn$, write $\nn_N \cq \nn_{1,N} + \nn_{2,N}$ and $\aa_N \cq \pp_N / (1 - \pp_N)$.
	Suppose that $\lim_N \nn_N \min\bra{ \pp_N, \: 1 - \pp_N } = \infty$.
	Then,
	for all $\cc \in \mbr$,
	we~have
	\[
		d_\TV^{\nn_{1,N}, \, \nn_{2,N}, \, \pp_N}\rbb{
			\rbb{ \tfrac12 \log(\aa_N \nn_N) + \cc } / \log\rbb{ \tfrac1{1 - \nn_{2,N} / \nn_N} }
		}
	\to
		2 \, \Phi(\tfrac12 e^{-\cc}) - 1
	\quad
		\asinf N.
	\]
\end{thm}

As in previous sections, for ease of presentation we omit the $N$-subscripts in the proof.
The technical calculations in this section are analogous to those in \S\ref{sec:rev:gibbs};
	the eigenfunctions are the same (after a reparametrisation)
but
	the eigenvalues are slightly different.

It is straightforward to check that the invariant distribution $m$ of the $\mcx$-chain is Binomial:
\[
	m(x)
=
	\binomt \nn\xx \pp^x (1 - \pp)^{\nn-\xx}
=
	\binomt \nn\xx \aa^\xx / (\aa+1)^\nn
\Qforall
	\xx \in \mcx.
\]
The eigenfunctions are then the family of polynomials orthogonal to the Binomial. These are the Krawtchouk polynomials (appropriately rescaled), defined precisely now.

\begin{defn}
\label{def:rev:gibbs:kraw}
	Define the \textit{Krawtchouk polynomials} $\bra{ K_i }_{i \in \mbn}$ via
	\[
		\Kab{\ii}{\xx}
	\cq
		\binom \nn\ii^{\!-1}
		\sumd[\min\bra{\ii,\xx}]{j=\max\bra{0,\ii-\nn-\xx}}
		\binom{\xx}{\jj} \binom{\nn-\xx}{\ii-\jj}
		\rbbb{-\frac1\aa}^\jj
	\Quad{for}
		\ii,\xx \in \mcx.
	\]
	When the second two parameters are fixed,
	abbreviate
	\[
		\varphi_\ii(\xx)
	\cq
		\Kab{\ii}{\xx}
	\Qfor
		\ii,\xx \in \mcx.
	\]
\end{defn}

The Krawtchouk polynomials are orthogonal with respect to the Binomial measure.


\begin{lem}[{\cite[\S1.10]{KS:hypergeo-askey}}]
\label{res:rev:gibbs:kraw-orthog}
	The Krawtchouk polynomials satisfy the orthogonality relations
	\[
		\sumt[\nn]{\xx=0}
		\Kab{\ii}{\xx} \Kab{\jj}{\xx} \aa^\xx \binomt \nn\xx
	=
		(\aa+1)^\nn \aa^{-\ii} \binomt N\ii^{-1} \delta_{\ii,\jj}
	\Qforall
		\ii,\jj \in \mcx.
	\]
	Thus the Krawtchouk polynomials are orthogonal with respect to the Binomial measure:
	\[
		\sumt[\nn]{\xx=0}
		m(\xx) \Kab{\ii}{\xx} \Kab{\jj}{\xx}
	=
		\aa^{-\ii} \binomt \nn\ii^{-1} \delta_{\ii,\jj}
	\Qforall
		\ii,\jj \in \mcx.
	\]
\end{lem}

The following proposition describes the eigenstatistics of this model; it is taken from \textcite[\S5.1]{DKSc:exp-families}.

\begin{prop}[Eigenstatistics; {\cite[\S5.1]{DKSc:exp-families}}]
\label{res:rev:gibbs:estats}
	The
		eigenvalues $\bra{\lambda_\ii}_{\ii \in \mcx}$
	and
		eigenfunctions\linebreak $\bra{f_\ii}_{\ii \in \mcx}$
	are given by the following:
	\[
		\lambda_\ii
	&\cq
		\prodd[\ii-1]{\jj=0}
		\frac{n_1 - \jj}{n_1 + n_2 - \jj}
	=
		\prodd[\ii-1]{\jj=0}
		\rbbb{ 1 - \frac{\nn_2}{\nn - \jj} }
	\ge
		0
	\Qfor
		\ii \in \mcx;
	\\
		f_\ii(\xx)
	&\cq
		\aa^{\ii/2} \binomt \nn\ii^{1/2} \Kab{\ii}{\xx}
	\equiv
		\aa^{\ii/2} \binomt \nn\ii^{1/2} \varphi_\ii(\xx)
	\Qfor
		\ii,\xx \in \mcx.
	\]
	Note that
		$f_\ii(0) = \aa^{\ii/2} \binomt \nn\ii^{1/2}$ for all $\ii \in \mcx$
	and
		$\lambda_i = 0$ for all $i \ge \nn_1 + 1$.
\end{prop}


Applying \cref{res:rev:gibbs:estats},
we obtain the following expressions for the terms in \cref{res:tv-approx:rev}:
\[
	\MT
&
\cq
	\sumt{\xx \in \mcx}
	m(\xx)
	\absb{ \sumt[M]{\ii=1} f_\ii(0) f_\ii(\xx) \lambda_\ii^\tt }
\\&
=
	(\aa+1)^{-\nn}
	\sumt[\nn]{\xx=0}
	\aa^\xx \binomt \nn\xx
	\absb{ \sumt[M]{\ii=1} \aa^\ii \binomt \nn\ii \varphi_\ii(\xx) \lambda_\ii^\tt };
\\
	\ET
&
\cq
	\sumt{\ii > \MM}
	\abs{ f_\ii(0) } \lambda_\ii^\tt
=
	\sumt{\ii > \MM}
	\alpha^{\ii/2} \binomt \nn\ii^{1/2} \lambda_\ii^\tt.
\]

Our first aim is to use this to determine which are the `important' eigenstatistics.

\begin{lem}[Error Term]
\label{res:rev:gibbs:error}
	For all $\eps > 0$ and all $\cc \in \mbr$,
	there exists an $\MM \cq \MM(\eps,\cc)$ so that,
	for $\tt \cq \tfrac12 \rbr{ \log(\aa \nn) + \cc } / \log\rbr{\tfrac1{1 - \nn_2/\nn}}$,
	if $I \cq \bra{1, ..., \MM}$, then
	\[
		\ET
	\le
		\ET'
	\le
		\eps
	\Qwhere
		\ET'
	\cq
		\sumt{\ii > \MM}
		\abs{ f_\ii(0) } \lambda_1^{\ii \tt}
	=
		\sumt{\ii > \MM}
		\alpha^{\ii/2} \binomt \nn\ii^{1/2}
		\prodt[\ii-1]{\jj=0} \rbb{ 1 - \tfrac{\nn_2}{\nn - \jj} }^\tt.
	\]
\end{lem}

\begin{Proof}
Observe that $0 \le \lambda_\ii \le \lambda_1^\ii = (1 - \nn_2/\nn)^\ii$ for all $\ii$.
The inequality $\ET \le \ET'$ now follows.
The equality in the definition of $\ET'$ is an immediate consequence of \cref{res:rev:gibbs:estats}.
For the inequality $\ET' \le \eps$,
choose $\MM$ so that $\sumt{\ii > \MM} e^{-\cc\ii} / \sqrt{\ii!} < \eps$;
then,
using \cref{res:rev:gibbs:estats} again,
we have
\[
	\sumt{\ii > \MM}
	\abs{ f_\ii(0) } \lambda_1^{\ii \tt}
\le
	\sumt{\ii > \MM}
	\rbb{ \sqrt{\aa \nn} \rbr{ 1 - \nn_2/\nn }^\tt }^\ii / \sqrt{\ii!}
=
	\sumt{\ii > \MM}
	e^{-\cc\ii} / \sqrt{\ii!}
<
	\eps.
\qedhere
\]
\end{Proof}

From now on, choose $\MM \cq \MM(\cc,\eps)$ as in \cref{res:rev:gibbs:error}.
Hence, for the main term, we need only deal with eigenstatistics with $i \asymp 1$.
We would then like to use the replacement $\lambda_\ii \approx (1 - \nn_2/\nn)^\ii$.

\begin{defn}[Adjusted Main Term]
\label{def:rev:gibbs:adj}
	Recalling that $\tt = \rbb{ \tfrac12 \log(\aa \nn) + \cc } / \log\rbb{\tfrac1{1 - \nn_2/\nn}}$,
	define
	\[
		\MT'
	\cq
		(\aa+1)^{-\nn} \sumt[\nn]{\xx=0} \aa^\xx \binomt \nn\xx
		\absb{ \sumt{\ii\ge1} \binomt \nn\ii \varphi_\ii(\xx) \aa^{\ii/2} e^{-\cc\ii} / n^{\ii/2} }.
	\]
\end{defn}

The following pair of lemmas approximate $\MT$ by $\MT'$ and then evaluate (asymptotically)~$\MT'$.

\begin{subtheorem}{thm}
	\label{res:rev:gibbs:main}

\begin{lem}[Main Term: Approximation]
\label{res:rev:gibbs:main:approx}
	For all $\eps > 0$ and $\cc \in \mbr$,
	for $\MM \cq \MM(\cc,\eps)$,
	we have
	\[
		\abs{ \MT - \MT' }
	\le
		2 \eps.
	\]
\end{lem}

It thus suffices to work with the $\MT'$, which has a significantly simpler form.
This is the main power of the technique: it allows us to replace the complicated $\lambda_i^\tt$ by the simpler $\lambda_1^{\ii \tt}$.
Typically, this power will be much easier to handle, particularly when melded with Binomial coefficients.

\begin{lem}[Main Term: Evaluation]
\label{res:rev:gibbs:main:eval}
	For all $\cc \in \mbr$,
	with $\MM \cq \MM(\cc,\eps)$,
	we have
	\[
		\tfrac12 \MT'
	\to
		2 \, \Phi\rbb{ \tfrac12 e^{-\cc} } - 1.
	\]
\end{lem}

\begin{Proof}[Proof of \cref{res:rev:gibbs:main:approx}]
Since $(1 - \nn_2/\nn)^\tt = (\aa \nn)^{1/2} e^{-\cc}$ and $\lambda_1 = 1 - \nn_2/\nn$,
we have
\[
	\abs{ \MT - \MT' }
&
=
	\absb{
		\sumt{\xx \in \mcx}
		m(\xx)
		\absb{
			\sumt[\MM]{\ii=1}
			f_\ii(0) f_\ii(\xx) \lambda_\ii^\tt
		}
	-	\sumt{\xx \in \mcx}
		m(\xx)
		\absb{
			\sumt[\infty]{\ii=1}
			f_\ii(0) f_\ii(x) \lambda_1^{\ii \tt}
		}
	}
\\&
\le
	\sumt{\xx \in \mcx}
	m(\xx)
	\sumt[\MM]{\ii=1}
	\abs{ f_\ii(0) f_\ii(\xx) } \abs{ \lambda_\ii^\tt - \lambda_1^{\ii \tt} }
+	\sumt{\xx \in \mcx}
	m(\xx)
	\sumt{\ii > \MM}
	\abs{ f_\ii(0) f_\ii(\xx) } \lambda_1^{\ii \tt}.
\]
We consider these two sums separately.
Recall that $\MM \cq \MM(\cc,\eps)$ is a constant.

For the first sum,
which we denote $S_1$,
we use the relation
\(
	\maxt{i \in [M]}
	\abs{ \lambda_\ii^\tt / \lambda_1^{\ii \tt} - 1 }
=
	\oh1,
\)
which is easy to derive.
Using
	Cauchy--Schwarz and the unit-normalisation of the eigenfunctions
as well as
	the relations
		$\lambda_1^{\ii \tt} = e^{-\cc \ii} (\aa \nn)^{-\ii/2}$
	and
		$f_\ii(0) = \alpha^{\ii/2} \binom \nn\ii^{1/2} \le (\alpha \nn)^{\ii/2} / \sqrt{\ii!}$,
we see that
\[
	S_1
&
=
	\sumt{\ii \in [\MM]}
	\abs{f_\ii(0)} \abs{ \lambda_\ii^\tt - \lambda_1^{\ii \tt} }
	\rbb{ \sumt{\xx \in \mcx} m(\xx) \abs{ f_\ii(\xx) } }
\\&
\le
	\maxt{\ii \in [\MM]} \abs{ \lambda_\ii^\tt / \lambda_1^{\ii \tt} - 1 }
\cdot
	\sumt{\ii \in [\MM]}
	\lambda_1^{\ii \tt} \abs{f_\ii(0)}
	\rbb{ \sumt{\xx \in \mcx} m(\xx) \abs{ f_\ii(\xx) } }
\\&
\le
	\oh1
\cdot
	\sumt{\ii \in [\MM]}
	e^{- \cc \ii} / \sqrt{\ii!}
=
	\oh1.
\]

For the second sum,
which we denote $S_2$,
using
	Cauchy--Schwarz and the unit-normalisation of the eigenfunctions again
and then
	the error term bound of \cref{res:rev:gibbs:error},
we see that
\[
	S_2
=
	\sumt{\ii > \MM}
	\abs{f_\ii(0)} \lambda_1^{\ii \tt}
	\rbb{ \sumt{\xx \in \mcx} m(\xx) \abs{f_\ii(\xx)} }
\le
	\sumt{\ii > \MM}
	\abs{f_\ii(0)} \lambda_1^{\ii \tt}
=
	\ET'
\le
	\eps.
\]

In conclusion, we see that
\(
	\abs{ \MT - \MT' }
\le
	\eps + \oh1
\le
	2 \eps
\)
(asymptotically),
as desired.
\end{Proof}

\begin{Proof}[Proof of \cref{res:rev:gibbs:main:eval}]
Evaluating this requires some algebraic manipulation then approximation.

For convenience, we drop some of the min/max from the limits in the sum in $\varphi_i$; define $\binom Nr \cq 0$ whenever it is not the case that $0 \le r \le N$.
Abbreviate $\zz \cq e^{-\cc}/\sqrt{\aa \nn}$.
For $\ell \in \bra{0, 1 , ..., \nn}$,
we~have
\[
	\sumd{i \ge 1}
	\binom \nn\ii \varphi_i(\ell) \frac{\aa^{i/2} e^{-\cc\ii}}{\nn^{i/2}}
&
=
	\sumd{i \ge 1}
	\sumd[\min\bra{i,\ell}]{\jj=0}
	\rbbb{-\frac1\aa}^\jj
	\binom \ell \jj \binom{\nn-\ell}{i-\jj} \rbr{ \aa \zz }^i
\\&
=
	\sumd{i \ge 0}
	\sumd[\min\bra{i,\ell}]{\jj=0}
	\rbbb{-\frac1\aa}^\jj
	\binom \ell \jj \binom{\nn-\ell}{i-\jj} \rbr{ \aa \zz }^i
-	1
\\&\textstyle
=
	\sumd[\ell]{\jj=0}
	\rbr{ -1/\aa }^\jj
	\binom \ell \jj \rbr{ \aa \zz }^\jj
	\sumd{i \ge \jj}
	\binom{\nn-\ell}{i-\jj} \rbr{ \aa \zz }^{i-\jj}
-	1
\\&
=
	\rbB{
		\sumt[\ell]{\jj=0}
		\rbr{ -1 }^\jj
		\binomt \ell \jj \zz^\jj
	}
\cdot
	\rbB{
		\sumt[\nn-\ell]{\rr=0}
		\binomt{\nn-\ell}{\rr} \rbr{ \aa \zz }^\rr
	}
-	1
\\&
=
	\rbr{ 1 - \zz }^\ell
\cdot
	\rbr{ 1 + \aa \zz }^{n-\ell}
-	1.
\]
We now need to take absolute values and average with respect to the weights $\aa^\ell \binomt \nn\ell / (\aa+1)^\nn$.

Observe that,
for any $\zeta \in \mbr$,
we have
\[
	\tfrac{\aa}{\aa+1}(1 - \zeta)
+	\tfrac{1}{\aa+1}(1 + \aa \zeta)
=
	1.
\]
So, setting $p_\xx \cq \tfrac{\aa}{\aa+1}(1 - \xx / \sqrt{\aa \nn})$ for $\xx \in \mbr$, the above is a $\Bin(\nn, p_{e^{-\cc}})$-type probability.
Indeed,
\[
	\MT'
&
=
	\sumt[\nn]{\ell=0}
	\absB{
		\binomt \nn\ell \cdot
			\rbb{ \tfrac{\aa}{\aa+1}(1 - \zz) }^\ell \cdot \rbb{ \tfrac{1}{\aa+1}(1 + \aa \zz) }^{\nn-\ell}
	-	\binomt \nn\ell \cdot
			\rbb{ \tfrac{\aa}{\aa+1} }^\ell \cdot \rbb{ \tfrac{1}{\aa+1} }^{\nn-\ell}
		}
\\&
=
	\sumt[\nn]{\ell=0}
	\absb{
		\binomt \nn\ell p_{e^{-\cc}}^\ell (1 - p_{e^{-\cc}})^{\nn-\ell}
	-	\binomt \nn\ell p_0^\ell (1 - p_0)^{\nn-\ell}
	}
=
	2 \, d_\TV\rbb{ \Bin\rbr{\nn, p_{e^{-\cc}}}, \: \Bin\rbr{\nn, p_0} }.
\]

It remains to compare these Binomials.
We do precisely this via the local CLT
in \cref{res:app:binom-tv}:
\[
	\tfrac12 \MT'
\to
	2 \, \Phi\rbb{ \tfrac12 e^{-\cc} } - 1
\quad
	\asinf \nn.
\qedhere
\]
\end{Proof}

\end{subtheorem}

We now have all the ingredients to establish the limit profile for this Gibbs sampler.

\begin{Proof}[Proof of \cref{res:rev:gibbs:res}]
Let us summarise what we have proved.
The following  are all evaluated at the target mixing time
\(
	\tt
=
	\tfrac12 \log(\aa \nn) + \cc \nn
\)
with $\MM \cq \MM(\cc,\eps)$ given by \cref{res:rev:gibbs:error}.
\begin{itemize}[noitemsep, topsep = \smallskipamount, label = \bcdot]
	\item 
	By \cref{res:rev:gibbs:error},
	we have
	\(
		\ET \le \eps.
	\)
	
	\item 
	By \cref{res:rev:gibbs:main:approx},
	we have
	\(
		\abs{\MT - \MT'} \le 2 \eps
	\)
	for $\nn$ sufficiently large
	
	\item 
	By \cref{res:rev:gibbs:main:eval},
	we have
	\(
		\tfrac12 \MT' \to 2 \, \Phi(\tfrac12 e^{-\cc}) - 1
	\)
	\asinf n.
\end{itemize}
Since $\eps > 0$ is arbitrary,
applying the TV-approximation lemma for reversible Markov chains, namely \cref{res:tv-approx:rev},
we immediately deduce the theorem.
\end{Proof}

\section{Random $k$-Cycle Walk on the Symmetric Group}
\label{sec:k-cycle}

\subsection{Walk Definition and Statement of Result}

We analyse the limit profile of the random $k$-cycle walk on the symmetric group $\symgr_n$.
This random walk starts (without loss of generality) from the identity permutation, and a step involves composing the current location with a uniformly chosen $k$-cycle.
This is an extension of the random transpositions studied by \textcite{T:limit-profile}.
We use representation theory for $k$-cycles, studied recently by \textcite{H:cutoff-k-cycle}, who established cutoff for any $2 \le k \ll n$, and found the order of the window if further $k \ll n/\log n$.
We determine the limit profile for any $2 \le k \ll n$.

For $\symgr_n$, the irreducible representations are indexed by partitions of $n$. As is common for card shuffles, the main contribution comes from those partitions with long first row; it is these we use as our set $I$.
We sharpen some of \citeauthor{H:cutoff-k-cycle}'s results slightly to determine the limit profile.

\begin{thm}[Random $k$-Cycle Walk]
\label{res:k-cycle:res}
	Let $n,k \in \mbn$.
	Consider the \textit{random $k$-cycle walk} on $\symgr_n$:
		start at $\id \in \symgr_n$;
		at each step, choose a $k$-cycle $\tau$ uniformly at random;
		move by right-multiplication.
	For $\tt \in \mbn_0$, write $d_\TV^{n,k}(\tt)$ for the TV distance of the random $k$-cycle walk on $\symgr_n$ from the uniform distribution on the appropriate set of permutations of a fixed parity, ie the odd ones if $k$ is even and $t$ is odd and the even ones otherwise.
	
	Let $(n_N)_\Ninn, (k_N)_\Ninn \in (\mbn \setminus \bra{1})^\mbn$.
	Suppose
		that $2 \le k_N \le n_N$ for all $\Ninn$
	and
		that $\lim_{N \to \infty} k_N / n_N = 0$.
	Then,
	for all $\cc \in \mbr$,
	we have
	\[
		d_\TV^{n_N, k_N}\rbb{ - n_N \rbr{ \log n_N + \cc } / \log(1 - k_N / n_N) }
	\to
		d_\TV\rbb{ \Pois(1 + e^{-\cc}), \: \Pois(1) }
	\quad
		\asinf N.
	\]
\end{thm}

Throughout the proof, for notational ease, we drop the subscripts, just writing $k$ and $n$, and assuming that $2 \le k \ll n$.
Write $\altgr_{n;k,t}$ for the set of odd permutations in $\symgr_n$ if $k$ is even and $t$ is odd and the even permutations otherwise.
Then, the $k$-cycle walk at time $t$ is supported on $\altgr_{n;k,t}$.

It is well-known that the irreducible representations for $\symgr_n$ are parametrised by partitions of $n$; see \cite{D:group-rep}.
We need to find a collection of irreducible representations which asymptotically contains all the total variation mass.
As is often the case with card shuffle-type walks, it is the partitions with long first row which we use.
More precisely,
for a partition $\lambda$ of $n$,
write $\lambda = (\lambda_1, ..., \lambda_n)$ with $\lambda_1 \ge \cdots \lambda_n$;
let $M \in \mbn$, and set
\[
	\mcp_n(M)
\cq
	\brb{ \lambda \text{ partition of } n \mid n - M < \lambda_1 < n }.
\]
The trivial representation, denoted $\triv^n$, corresponds to the partition with only one block, ie $\triv^n_1 = n$ and $\triv^n_i = 0$ for $i \ge 2$.
Write $\mcp^*_n(M) \cq \mcp_n(M) \cup \bra{\triv^n}$.

We now phrase \citeauthor{T:limit-profile}'s lemma, ie \cref{res:tv-approx:tey}, in this set-up.

\begin{lem}
\label{res:k-cycle:tey}
	For all $t \in \mbn_0$ and all $M \ge 1$,
	we have
	\[
		\absb{
			d_\TV(\tt)
		-	\tfrac1{n!} \sumt{\sigma \in \altgr_{n;k,t}} \absb{ \sumt{\lambda \in \mcp_n(M)} d_\lambda s_\lambda(k)^t \chi_\lambda(\sigma) }
		}
	\le
		\sumt{\lambda \notin \mcp^*_n(M)} d_\lambda \abs{s_\lambda(k)}^t
	\label{eq:k-cycle:tv-mt-et}
	\]
\end{lem}

Given $k$ and $t$, the random walk is supported on the set of permutations with a fixed sign; half the permutations are odd and half are even.
Hence the factor $\tfrac12 \abs{\altgr_{n;k,t}}^{-1} = \tfrac1{n!}$ in the lemma above.
(We emphasise the dependence on $k$ in the character ratio $s_\lambda(k)$.)

\begin{Proof}[Outline of Proof of \cref{res:k-cycle:res}]
We show in \cref{res:k-cycle:error} that,
for
	all $\eps > 0$
and
	all $\cc \in \mbr$,
there exists a constant $M \cq M(\eps,\cc)$ so that
this the right hand side of \eqref{eq:k-cycle:tv-mt-et} is at most $\eps$ when $t = - n \rbr{ \log n + \cc } / \log(1 - k/n)$.
Thus, for the main term, we are interested in $\lambda$ with $n - \lambda_1 \asymp 1$.

It is well-known that $d_\lambda \le \binomt n{\lambda_1} d_{\lambda^*}$, where $\lambda^* \cq \lambda \setminus \lambda_1$ is the partition $\lambda$ with the largest element removed.
In fact, $d_\lambda \approx \binomt nr d_{\lambda^*} \approx \tfrac1{r!} n^r d_{\lambda^*}$ when $r \cq n - \lambda_1 \asymp 1$; see \cite[Proposition~3.2]{T:limit-profile}.

\textcite[Theorem~5]{H:cutoff-k-cycle} states a rather general result on the character ratios $s_\lambda(k)$.
Manipulating this general formula in the special case of $\lambda \in \mcp_n(M)$, ie $r = n - \lambda_1 \asymp 1$,
we show in \cref{res:k-cycle:long} that
\(
	s_\lambda(k)
\approx
	(1 - k/n)^r.
\)
When raised to the power $t$, we get
\(
	s_\lambda(k)
\approx
	n^{-r} e^{-r \cc}.
\)

Altogether, by allowing us to replace
\[
	\sumt{\lambda \in \mcp_n(M)}
	d_\lambda s_\lambda(k)^t \chi_\lambda(\sigma)
\Qwith
	\sumt[M-1]{r=1}
	\tfrac1{r!} e^{-r \cc}
	\sumt{\lambda : \lambda_1 = n - r}
	d_{\lambda^*} \chi_\lambda(\sigma),
\]
this converts an unmanageable main term sum into what is in essence a generating function.
We then adapt results of \textcite[\S4]{T:limit-profile} to control this generating function.
\end{Proof}

As stated above, to prove this theorem we use representation theory results on the $k$-cycle walk from \cite{H:cutoff-k-cycle}.
We state these precisely in the next section; we have to sharpen some results slightly.
Throughout this section, $\lambda$ will always be a partition of $n$, written $\lambda \vdash n$.

\medskip

Following \cite{H:cutoff-k-cycle},
we use the \textit{Frobenius notation} for a partition:
\[
	\lambda = (a_1, ..., a_m \mid b_1, ..., b_m)
\Qwith
	a_i \cq \lambda_i - i + \tfrac12
\Qand
	b_i \cq \lambda'_i - i + \tfrac12,
\]
where $\lambda'$ is the transpose of the partition $\lambda$.
Writing $r \cq n - \lambda_1$,
the following hold:
\begin{gather*}
	a_1 = n - r - 1,
\quad
	a_1 - a_i = n - (1 + r + \lambda_i - i),
\\
	a_1 + b_i = n - (r - \lambda'_i + i)
\Quad{and}
	\max\bra{a_2, ..., a_m, b_1, b_2, ..., b_m} \le r.
\end{gather*}
We use the following notation for the \textit{descending factorial}:
for $z \in \mbr$ and $k \in \mbn$,
write
\[
	\defac zk
\cq
	z (z - 1) \cdots (z - k + 1).
\]

Without further ado, we quote the required results from \textcite{H:cutoff-k-cycle} in the next subsection.

\subsection{Statements of Character Ratio Bounds}
\label{sec:k-cycle:statements}

In this subsection, we state a result from \cite{H:cutoff-k-cycle}, and deduce some corollaries of these statements.
We do not give any proofs at this stage; these are deferred to \S\ref{sec:k-cycle:proofs}.

\medskip

The first result which we quote determines asymptotically the character ratio for partitions with long first row---which, we recall, are the partitions of particular interest to us.

\begin{thm}[{\cite[Theorem~5(a)]{H:cutoff-k-cycle}}]
\label{res:k-cycle:5a}
	Let $0 < \eps < \tfrac12$.
	Suppose that $r + k + 1 < \tfrac13 n$.
	Then
	\[
		s_\lambda(k)
	&
	=
		\frac{ \defac{n-r-1}{k} }{ \defac nk }
	\cdot
		\prod_{i=2}^m \rbbb{ 1 - \frac{k}{n - (1 + r + \lambda_i - i)} }
	\cdot
		\prod_{i=1}^m \rbbb{ 1 - \frac{k}{n - (r - \lambda'_i + i)} }^{\!-1}
	\\&\qquad
	+	\Ohb{ \expb{ k \log(k+r+1) + \Oh{1/\sqrt r} - k \log(n-k) } }.
	\]
	Further, if $r < k$, then the error term is actually 0.
\end{thm}

In this article, we are interested in partitions with long first row, namely $\mcp_n(M)$.
We can apply this theorem to analyse asymptotics of partitions with long first row.
We defer the proof to \S\ref{sec:k-cycle:proofs}.

\begin{cor}[Long First Row]
\label{res:k-cycle:long}
	Let $2 \le k \le \tfrac13 n$.
	Let $r \in \mbn$.
	Let $\lambda \vdash n$ with $\lambda_1 = n - r$.
	Then,
	\[
		s_\lambda(k)
	=
		\rbr{ 1 - k/n }^r \cdot \rbb{ 1 + \Oh{k/n^2} }.
	\]
\end{cor}

This covers the case where the first row is long.
The next two results consider shorter rows; the first is for $k \ge 6 \log n$ and the second for $k \le 6 \log n$.
These statements are not exactly the same as in \cite{H:cutoff-k-cycle}, but are slight strengthenings; their proofs are given in \S\ref{sec:k-cycle:proofs}.

\begin{thm}[cf {\cite[Theorem~5(b)]{H:cutoff-k-cycle}}]
\label{res:k-cycle:short:k-big}
	Assume that $6 \log n \le k \ll n$.
	Let $\theta \cq 0.68 > \tfrac23$; so $e^{-\theta} > 0.506$.
	Consider $\lambda$ with $b_1 \le a_1 \le e^{-\theta} n$.
	Then,
	\[
		\abs{s_\lambda(k)} \le \expb{ - \rbr{\tfrac12 + \tfrac1{10}} k }.
	\]
\end{thm}

\begin{lem}[cf {\cite[Lemmas~14 and 15]{H:cutoff-k-cycle}}]
\label{res:k-cycle:short:k-small}
	Assume that $2 \le k \le 6 \log n$.
	Let $\lambda \vdash n$ with $b_1 \le a_1$ and $r \cq n - \lambda_1 \in [\tfrac13 n, n]$.
	Then,
	\[
		\abs{s_\lambda(k)}
	\le
		\expb{ - (\tfrac12 + \tfrac1{10}) rk/n}.
	\]
\end{lem}

From these statements, along with the standard bounds on $d_\rho$, the dimension of an irreducible representation $\rho$, we are able to control the two terms, which we call the \textit{main} and \textit{error} terms, in \cref{res:tv-approx:tey}.
Our first port of call is to find a suitable $M$ to bound the error term.
Once we have determined this, for the main term we need only consider partitions $\lambda$ with $\lambda_1 \ge n - M$.
We take $M$ to be order 1 (but arbitrarily large); so $\lambda_1 \ge n - M$ falls into the ``long first row'' case.

\begin{lem}[Error Term]
\label{res:k-cycle:error}
	Let $\cc \in \mbr$ and $t \cq - n(\log n + \cc) / \log(1 - k/n)$.
	For $M \in \mbn$, let
	\[
		\ET_M
	\cq
		\sumt{\lambda : \lambda_1 \le n - M}
		d_\lambda \abs{s_\lambda(k)}^t
	=
		\sumt{r \ge M}
		\sumt{\lambda : \lambda_1 = n - r}
		d_\lambda \abs{s_\lambda(k)}^t.
	\]
	Then,
	\(
		\ET_M \to 0 \ \asinf M.
	\)
\end{lem}

This controls the error term.
We now consider the main term in \cref{res:k-cycle:tey}.

\begin{lem}[Main Term]
\label{res:k-cycle:main}
	Let $\cc \in \mbr$ and $t \cq - n(\log n + \cc) / \log(1 - k/n)$.
	For $M \in \mbn$, let
	\[
		\MT_M
	\cq
		\tfrac1{n!} \sumt{\sigma \in \altgr_{n;k,t}} \absb{ \sumt{\lambda \in \mcp_n(M)} d_\lambda s_\lambda^t \chi_\lambda(\sigma) }.
	\]
	Then,
	\(
		\MT_M \to d_\TV\rbb{ \Pois(1 + e^{-\cc}), \: \Pois(1) } \ \asinf M.
	\)
\end{lem}

\begin{Proof}[Proof of \cref{res:k-cycle:res} Given \cref{res:k-cycle:error,res:k-cycle:main}]
\cref{res:k-cycle:tey} formulates \citeauthor{T:limit-profile}'s lemma, ie \cref{res:tv-approx:tey}, in the set-up of the random $k$-cycle walk.
\cref{res:k-cycle:error,res:k-cycle:main} control the error and main terms, respectively.
Combining these three ingredients establishes \cref{res:k-cycle:res}.
\end{Proof}

It remains to control the error and main terms, ie prove \cref{res:k-cycle:error,res:k-cycle:main} respectively.

\subsection{Controlling the Main and Error Terms}
\label{sec:k-cycle:mt-et}

We control the main term in \S\ref{sec:k-cycle:mt-et:mt} and the error term in \S\ref{sec:k-cycle:mt-et:et}.

\subsubsection{Controlling the Main Term}
\label{sec:k-cycle:mt-et:mt}

We analyse the main term, ie \cref{res:k-cycle:main}, first.
The analysis follows similarly to the case of random transpositions (ie $k = 2$) considered by \textcite{T:limit-profile}.
We need only consider partitions $\lambda$ with long first row, namely $\lambda_1 = n - r$ with $1 \le r \le M$, where $M$ is some (arbitrarily large) constant.
These are precisely the partitions considered in the results quoted from \textcite{H:cutoff-k-cycle}.

\textcite[\S4.1 and \S4.2]{T:limit-profile} then has some technical lemmas to get the main term into the desired form.
We summarise these now.
Note that he considers time $\tfrac12 n \log n + \cc n = \tfrac12 n (\log n + 2 \cc)$, while we are considering $t = - n \rbr{ \log n + \cc } / \log(1 - k/n)$; hence our two $\cc$-s differ by a factor 2.

Before digging into the details of his lemmas, we give the high-level reasons why his proof passes over to our case.
When considering the main term, one need only study those partitions with long first row, ie $\lambda$ with $r \cq n - \lambda_1 \asymp 1$.
For such $\lambda$, consider the difference between $s_\lambda(2)$ and $s_\lambda(k)$:
\[
	s_\lambda(2)
&
=
	\rbr{ 1 - 2/n }^k \cdot \rbb{ 1 + \Oh{1/n^2} }
=
	\expb{-2r/n} \cdot \rbb{ 1 + \Oh{1/n^2} }
\\
	s_\lambda(k)
&
=
	\rbr{ 1 - k/n }^r \cdot \rbb{ 1 + \Oh{k/n^2} }
=
	\expb{-kr/n} \cdot \rbb{ 1 + \Oh{k^2/n^2} }.
\]
\citeauthor{T:limit-profile} needs $s_\lambda(2)^\tt = n^{-r} e^{-\cc}$.
This goes some way to justifying why we \emph{expect} $\tt \approx \tfrac nk \log n$ to be the mixing time, and that the cutoff window should scale down with $k$ linearly.

\medskip

We now proceed more formally.
For each $r \ge 1$,
define the polynomials $T_r$ by
\[
	T_r(z)
\cq
	\sumd[r]{\ii=1}
	\binom z{r-\ii} \frac{(-1)^\ii}{\ii!}
\Qfor
	z \in \mbn_0.
\]
For a partition $\lambda$,
write $\lambda^* \cq \lambda \setminus \lambda_1$ for $\lambda$ with the first row removed.
For a permutation $\sigma \in \symgr_n$, write $\Fix \sigma$ for the number of fixed points in $\sigma$.

\begin{lem}[{\cite[Lemma~4.3]{T:limit-profile}}]
\label{res:k-cycle:4.3}
	Let $r \in \mbn$. Let $\sigma \in \symgr_n$ be a permutation with at least one cycle of length greater than $r$.
	Then,
	\[
		\tfrac1{r!}
		\sumt{\lambda \vdash n : \lambda_1 = n - r}
		d_{\lambda^*} \chi_\lambda(\sigma)
	=
		T_r(\Fix \sigma).
	\]
\end{lem}

The proof of this lemma is combinatorial and strongly relies on the Murnaghan--Nakayama rule.
\cref{res:k-cycle:4.3} is a statement purely about the representation theory of the symmetric group; it is nothing to do with the random walk.
%
Using this result, one can obtain the following approximation.

\begin{lem}[cf {\cite[Lemma~4.2]{T:limit-profile}}]
\label{res:k-cycle:4.2}
	Set $t \cq - n \rbr{ \log n + \cc } / \log(1 - k/n)$.
	Let $M \in \mbn$.~%
	Then,
	\[
		\frac1{n!}
		\sumd{\sigma \in \altgr_{n;k,t}}
		\absbb{
			\sumd[M-1]{r=1}
			\sumd{\lambda : \lambda_1 = n - r}
			d_\lambda s_\lambda(k)^t \chi_\lambda(\sigma)
		}
	=
		\frac1{n!}
		\sumd{\sigma \in \altgr_{n;k,t}}
		\absbb{
			\sumd[M-1]{r=1}
			e^{-r\cc} T_r(\Fix \sigma)
		}
	+	\oh1.
	\]
\end{lem}


	%
To prove this, one separates $\altgr_{n;k,t}$ into the set of permutations with a cycle of length greater than $M$ and those with all cycles of length at most $M$.
Also, it is not difficult to check, using the hook-length formula and \cref{res:k-cycle:long} (cf \cite[Propositions~3.1 and 3.2]{T:limit-profile}), that
\[
	\absb{ d_\lambda s_\lambda(k)^t - e^{-r\cc} d_{\lambda^*} / r! }
=
	\Oh{ \log n / n }
\Quad{when}
	t = - n (\log n + \cc) / \log(1 - k/n).
\]
This is the crucial bound in \cite[Lemma~4.2]{T:limit-profile}.
The remainder of that proof uses only facts about the symmetric group, not specific to $k$-cycles, or the time $t$.
Hence, the result follows through.
Next, polynomials of high degree can be neglected polynomials, in the following sense.

\begin{lem}[{\cite[Lemma~4.4]{T:limit-profile}}]
\label{res:k-cycle:4.4}
	For any $M$ with \toinf M \asinf n, we have
	\[
		\absbb{
			\frac1{n!}
			\sumd{\sigma \in \altgr_{n;k,t}}
			\absbb{
				\sumd[M-1]{r=1}
				e^{-r\cc} T_r(\Fix \sigma)
			}
		-	\frac1{n!}
			\sumd{\sigma \in \altgr_{n;k,t}}
			\absbb{
				\sumd[\infty]{r=1}
				e^{-r\cc} T_r(\Fix \sigma)
			}
		}
	\to
		0.
	\]
\end{lem}

We must next evaluate this infinite sum.
For $\cc \in \mbr$, define the function $f_\cc$ by
\[
	f_\cc
:
	m \mapsto \expb{-e^{-\cc}} \rbb{1 + e^{-\cc}}^m - 1
:
	\mbn \to \mbr.
\]

\begin{prop}[{\cite[Proposition~4.5]{T:limit-profile}}]
\label{res:k-cycle:4.5}
	Let $m \in \mbn$.
	Then,
	\[
		\sumt[\infty]{r=1}
		e^{-r\cc} T_r(m)
	=
		f_\cc(m).
	\]
\end{prop}

Finally, we evaluate this function at $\Fix \sigma$ with $\sigma \sim \Unif(\altgr_{n;k,t})$ and take the expectation.

\begin{lem}[{cf \cite[Lemma~4.6]{T:limit-profile}}]
\label{res:k-cycle:4.6}
	We have
	\[
		\tfrac2{n!}
		\sumt{\sigma \in \altgr_{n;k,t}}
		\abs{ f_\cc(\Fix \sigma) }
	\to
		\exb{ \abs{ f_\cc(\Pois(1)) } }
	=
		2 \, d_\TV\rbb{ \Pois(1 + e^{-\cc}), \: \Pois(1) }.
	\]
\end{lem}

The idea behind this lemma is simple:
	it is well-known that if $\sigma \sim \Unif(\symgr_n)$ then $\Fix \sigma \to^d \Pois(1)$;
	we show that the same is true when $\sigma$ is restricted to having a prescribed parity.

\begin{Proof}[Proof of \cref{res:k-cycle:4.6}]
We claim that precisely half the permutations with a given number of fixed points are even (and hence half are odd):
if $\altgr_n$ is the alternating group of even permutations,
then
\[
	\absb{ \brb{ \sigma \in \altgr_n \subseteq \symgr_n \mid \Fix \sigma = r } }
=
	\tfrac12 \absb{ \brb{ \sigma \in \symgr_n \mid \Fix \sigma = r } }
\Qforall
	r \ge 0.
\]
Given this claim, the lemma follows easily, as in \cite[Lemma~4.6]{T:limit-profile}.

We now justify our claim.
First, we find the number of permutations in $\symgr_n$ (of either parity) with \emph{exactly} $r$ fixed points, which we denote $f_{n,r}$.
Note that
\(
	f_{n,r} = \binom nr f_{n-r,0}.
\)
Indeed:
	first select the $r$ points to be fixed, for which there are $\binom nr$ choices;
	then choose a permutation on the remaining $n-r$ points with no fixed points.
It remains to calculate $f_{m,0}$, ie the number of \textit{derangements} of $m$ objects, for each $m \in \bra{0, ..., n}$.
To do this, we use the inclusion--exclusion principle.
For $i \in [m]$, let
\(
	\symgr_{m,i}
\cq
	\bra{ \sigma \in \symgr_m \mid \sigma(i) = i }
\)
denote the set of permutations on $m$ objects that fix the $i$-th object.
Observe that
\(
	\abs{\cap_{i \in I} \symgr_{m,i}}
=
	\abs{\symgr_{m - \abs I}}
=
	(m - \abs I)!
\)
for all $I \subseteq [m]$; for each $\ell \in [m]$, there are $\binom m\ell$ choices of $I \subseteq [m]$ with $\abs I = \ell$.
Hence, by inclusion--exclusion, we have
\[
	f_{m,0}
=
	m! \sumt[m]{\ell=0} (-1)^\ell / \ell!.
\]
Combined with the fact that $f_{n,r} = \binom nr f_{n-r,0}$, we thus deduce that
\[
	f_{n,r}
=
	\binomt nr \cdot (n-r)! \sumt[n-r]{\ell=0} (-1)^\ell / \ell!
=
	\tfrac{n!}{r!} \sumt[n-r]{\ell=0} (-1)^\ell / \ell!.
\]

We now turn to even permutations, ie $\altgr_n$.
We apply an analogous method.
Denote by $f'_{n,r}$ the number of permutations in $\altgr_n$ with \emph{exactly} $r$ fixed points.
Since appending fixed points to a permutation does not change its parity, again we have
\(
	f'_{n,r} = \binom nr f'_{n-r,0}.
\)
For $m \in \mbn$ and $i \in [m]$, define
\(
	\altgr_{m,i}
\cq
	\bra{ \sigma \in \altgr_m \mid \sigma(i) = i }.
\)
Analogously to before,
since appending fixed points does not affect the parity,
we have
\(
	\abs{\cap_{i \in I} \altgr_{m,i}}
=
	\abs{\altgr_{m - \abs I}}
=
	\tfrac12 (m - \abs I)!
\)
for all $I \subseteq [m]$.
This~is~a~factor~$\tfrac12$ different to $\abs{\cap_{i \in I} \symgr_{m,i}}$ from before.
Using the inclusion--exclusion principle thus gives, as~before,
\[
	f'_{m,0}
=
	\tfrac12 \sumt[m]{\ell=0} (-1)^\ell / \ell!
\Qand
	f'_{n,r} = \tfrac12 \tfrac{n!}{r!} \sumt[n-r]{\ell=0} (-1)^\ell / \ell!.
\]
Since half the permutations of a given sign are even, ie $f'_{n,r} = \tfrac12 f_{n,r}$, the other half must be odd.
\end{Proof}

Observe that \cref{res:k-cycle:4.4,res:k-cycle:4.5,res:k-cycle:4.6} are statements purely about the representation theory of the symmetric group; it is nothing to do with the random walk.

Using standard applications of the triangle inequality, these lemmas can then be combined to deduce that the main term converges to the TV-distance in question; see \cite[\S4.4]{T:limit-profile}.

\begin{Proof}[Proof of \cref{res:k-cycle:main}]
Let $\eps > 0$ and let $M$ and $n$ be large enough so that all the approximations are true up to an additive error of $\eps$.
The following inequalities hold:
\begin{alignat*}{2}
	&\makebox[7em][l]{by \cref{res:k-cycle:tey,res:k-cycle:error},}
&\quad
	\absbb{
		d_\TV(\tt)
	-	\frac1{n!} \sumd{\sigma \in \altgr_{n;k,t}} \absB{ \sumd{\lambda \in \mcp_n(M)} d_\lambda s_\lambda(k)^t \chi_\lambda(\sigma) }
	}
&\le
	\eps;
\\
	&\makebox[7em][l]{by \cref{res:k-cycle:4.2},}
&\quad
	\absbb{
		\frac1{n!} \sumd{\sigma \in \altgr_{n;k,t}} \absB{ \sumd{\lambda \in \mcp_n(M)} d_\lambda s_\lambda(k)^t \chi_\lambda(\sigma) }
	-	\frac1{n!} \sumd{\sigma \in \altgr_{n;k,t}} \absB{ \sumd[M-1]{r=1} e^{-r\cc} T_r(\Fix \sigma) }
	}
&\le
	\eps;
\\
	&\makebox[7em][l]{by \cref{res:k-cycle:4.4},}
&\quad
	\absbb{
		\frac1{n!}
		\sumd{\sigma \in \altgr_{n;k,t}}
		\absB{
			\sumd[M-1]{r=1}
			e^{-r\cc} T_r(\Fix \sigma)
		}
	-	\frac1{n!}
		\sumd{\sigma \in \altgr_{n;k,t}}
		\absB{
			\sumd[\infty]{r=1}
			e^{-r\cc} T_r(\Fix \sigma)
		}
	}
&\le
	\eps;
\\
	&\makebox[7em][l]{by \cref{res:k-cycle:4.6},}
&\quad
	\absbb{
		\frac1{n!}
		\sumd{\sigma \in \altgr_{n;k,t}}
		\abs{ f_\cc(\Fix \sigma) }
	-	d_\TV\rbb{ \Pois(1 + e^{-\cc}), \: \Pois(1) }
	}
&\le
	\eps.
\end{alignat*}
Since $\eps > 0$ was arbitrary, the proof is now complete by the triangle inequality.
\end{Proof}

\subsubsection{Controlling the Error Term}
\label{sec:k-cycle:mt-et:et}

Finally, we control the error term, ie \cref{res:k-cycle:error}.
As before, we consider only $\lambda$ with $\lambda_1 \ge \lambda'_1$.
%
Consider first the dimensions of the irreducible representations, ie $d_\lambda$.

\begin{lem}
\label{res:k-cycle:et:irreps}
	The following bounds hold:
	\begin{alignat}{3}
		\sumt{\lambda \vdash n : \lambda_1 = n - r} d_\lambda
	&\le
		n^r 2^r / r^{r/2}
	&&\Quad{for}&
		r &\in [1, n];
	\label{eq:k-cycle:irreps:nr/rr}
	\\
		\sumt{\lambda \vdash n : \lambda_1 = n - r} d_\lambda
	&\le
		n^{r/2} 4^r
	&&\Quad{for}&
		r &\in [\tfrac13 n, n].
	\label{eq:k-cycle:irreps:nr4r}
	\end{alignat}
\end{lem}

\begin{Proof}
It is well-known that $d_\lambda$ is equal to the number of ways of placing the numbers $1$ through $n$ into the \textit{Young diagram} of $\lambda$; see, eg, \cite[Lemma~6]{DS:random-trans}.
From this, it is immediate that $d_\lambda \le \binom n{\lambda_1} d_{\lambda^*}$ where $\lambda^* \cq \lambda \setminus \lambda_1$ is the partition obtained by removing the largest element of $\lambda$.
It is also standard that $\sum_{\rho \vdash r} d_\rho^2 = \abs{\mcs_r} = r!$; see, eg, \cite[Theorem~3.8.11]{CsST:harmonic-analysis-finite-groups}.
(This last claim is true for any group, not just the symmetric group.)
Associate to the partition $\rho = (\rho_1, ..., \rho_r)$ of $[r]$, written in increasing order, the subset $\tilde \rho \cq \bra{\rho_1, \rho_1 + \rho_2, ..., \rho_1 + \cdots + \rho_r}$ of $[r]$.
This mapping is injective, and so
\(
	\abs{\bra{\rho \mid \rho \vdash r}}
\le
	\abs{\bra{\tilde \rho \mid \tilde \rho \subseteq [r]}}
=
	2^r.
\)
Combining these bounds and using Cauchy--Schwarz~gives
\[
	\sumt{\lambda \vdash n: \lambda_1 = n-r}
	d_\lambda
\le
	\binomt nr
	\sumt{\rho : \rho \vdash r}
	d_\rho
\le
	\binomt nr
	\sqrt{ \sumt{\rho \vdash r} 1 \cdot \sumt{\rho \vdash r} d_\rho^2 }
\le
	\defac nr \rbr{ 2^r r! }^{1/2} / r!
\le
	n^r 2^r / r^{r/2}.
\]

The second claim is a special case:
\(
	n^r / r^{r/2}
\le
	n^r / (\tfrac13 n)^{r/2}
=
	3^{r/2} n^{r/2}
\le
	2^r n^{r/2}
\)
when
\(
	r \ge \tfrac13 n.
\)
\end{Proof}

We split the summation $\sumt[\infty]{r=M}$ in the error term $\ET_M$ into two parts:
	$r \le 0.495 n$ and $r \ge 0.495 n$;
	the latter sum is separated according to whether or not $k \le 6 \log n$.

\begin{Proof}[Proof of \cref{res:k-cycle:error}]
Throughout this proof,
let $\cc \in \mbr$ and
\(
	t \cq - n \rbr{ \log n + \cc } / \log(1 - k/n).
\)

\emph{Consider first $r \in [M,0.495n]$ with $2 \le k \ll n / \log n$.}
Recall \cref{res:k-cycle:long}
which implies that
\[
	\abs{s_\lambda(k)}^t
\le
	\expb{ t r \log(1 - k/n) + \Oh{ t k/n^2 } }
\Qwhere
	r \cq n - \lambda_1.
\]
Note that $t k \asymp n \log n \ll n^2$.
Thus,
for all $\cc \in \mbr$,
using \eqref{eq:k-cycle:irreps:nr/rr},
we have
\begin{subequations}
	\label{eq:k-cycle:et:E1}
\begin{align}
	\mce_1
&
\cq
	\sumt[0.495n]{r=M}
	\sumt{\lambda_1 = n-r}
	d_\lambda \abs{s_\lambda(k)}^t
\\&\phantom{:}
\le
	\sumt[\infty]{r=M}
	n^r 2^r r^{-r/2} \cdot 2 \expb{ - n (\log n + \cc) \cdot r/n }
\\&\phantom{:}
\le
	\sumt[\infty]{r=M}
	2 \expb{ \tfrac12 r \rbb{ \abs \cc - \log(r/4) } }.
\end{align}
\end{subequations}
The summand is independent of $n$, and gives rise to a summable series; hence $\mce_1 \to 0$ as \toinf M.

\emph{Consider next $r \in [0.495n, n]$ with $6 \log n \le k \ll n / \log n$.}
When $r \ge 0.495n$ we have $a_1 \le 0.505n < e^{-0.68}n$.
Recall \cref{res:k-cycle:short:k-big}
which implies that
\[
	\abs{s_\lambda(k)}^t
\le
	\expb{-(\tfrac12 + \tfrac1{10})kt}.
\]
Now, $tk = n(\log n + \cc)$.
Hence,
for all $\cc \in \mbr$,
using \eqref{eq:k-cycle:irreps:nr4r},
we have
\begin{subequations}
	\label{eq:k-cycle:et:E2>}
\begin{align}
	\mce_{2,>}
&
\cq
	\sumt[n]{r=0.495n}
	\sumt{\lambda_1 = n-r}
	d_\lambda \abs{s_\lambda(k)}^t
\\&\phantom{:}
\le
	\sumt[n]{r=0.495n}
	4^r n^{r/2} \expb{ -(\tfrac12 + \tfrac1{10}) \cdot n(\log n + \cc) }
\\&\phantom{:}
\le
	\sumt[n]{r=1}
	n^{-1} \expb{ \tfrac12 r \log n - \tfrac12 n \log n - \tfrac1{10} n \log n + n (\abs \cc + 1) }
\\&\phantom{:}
\le
	\expb{ n(\abs \cc + 1) - \tfrac1{10} n \log n }
=
	\oh1.
\end{align}
\end{subequations}

\emph{Consider finally $r \in [0.495n, n]$ with $2 \le k \le 6 \log n$.}
Note that $r \ge 0.495 n \ge \tfrac13 n$.
Recall \cref{res:k-cycle:short:k-small}
which implies that
\[
	\abs{s_\lambda(k)}^t
\le
	\expb{-(\tfrac12 + \tfrac1{10})tkr/n}.
\]
Now, $t \ge \tfrac{9}{10} \tfrac nk \log n$, as $-1/\log(1-x) > \tfrac{9}{10} x^{-1}$ for $x < \tfrac1{10}$.
Hence,
for all $\cc \in \mbr$,
using \eqref{eq:k-cycle:irreps:nr4r},
we~have%
\begin{subequations}
	\label{eq:k-cycle:et:E2<}
\begin{align}
	\mce_{2,<}
&
\cq
	\sumt[n]{r=0.495n}
	\sumt{\lambda_1 = n-r}
	d_\lambda \abs{s_\lambda(k)}^t
\\&\phantom{:}
\le
	\sumt[n]{r=0.495n}
	4^r n^{r/2} \expb{ -r \log n \cdot \tfrac{9}{10}(\tfrac12 + \tfrac1{10}) }
\\&\phantom{:}
\le
	\sumt[n]{r=0.495n}
	n^{-1} \expb{ \tfrac12 r \log n - \tfrac12 r \log n - \tfrac1{25} r \log n + n \log 4 }
\\&\phantom{:}
\le
	\expb{ n \log 4 - \tfrac1{60} n \log n }
=
	\oh1.
\end{align}
\end{subequations}

The lemma follows immediately from these three considerations, namely (\ref{eq:k-cycle:et:E1}, \ref{eq:k-cycle:et:E2>}, \ref{eq:k-cycle:et:E2<}).
Indeed, define $\mce_2 \cq \mce_{2,<} \one{2 \le k \le 6 \log n} + \mce_{2,>} \one{6 \log n < k \ll n / \log n}$.
Then
\(
	\ET_M = \mce_1 + \mce_2.
\)
The lemma follows from (\ref{eq:k-cycle:et:E1}, \ref{eq:k-cycle:et:E2>}, \ref{eq:k-cycle:et:E2<}) combined.
\end{Proof}

\subsection{Proofs of Character Ratio Bounds}
\label{sec:k-cycle:proofs}

In this section we give the deferred proofs from \S\ref{sec:k-cycle:statements}.

\begin{Proof}[Proof of \cref{res:k-cycle:long}]
Write $P_0$, $P_1$ and $P_2$ for the three terms in the product from \cref{res:k-cycle:5a}:%
\[
	P_0
\cq
	\frac{\defac{n-r-1}{k}}{\defac nk};
\quad
	P_1
\cq
	\prod_{i=2}^m \rbbb{ 1 - \frac{k}{n - (1 + r + \lambda_i - i)} };
\quad
	P_2
\cq
	\prod_{i=1}^m \rbbb{ 1 - \frac{k}{n - (r - \lambda'_i + i)} }^{\!-1}.
\]
Then, by \cref{res:k-cycle:5a} (ie \cite[Theorem~5(a)]{H:cutoff-k-cycle}), the main contribution to $s_\lambda(k)$ is $P_0 P_1 P_2$.

Since $r$ is a constant, all the $\bra{a_i, b_i, \lambda_i, \lambda'_i}$ are order 1, with the exception of $\lambda_1 = n - r$ and $a_1 = n - r - \tfrac32$.
Hence, all the terms in the two products are very similar to $1-k/n$ and in the first term to $1-r/n$.
In particular, for $1 \le j,\ell \le \tfrac12 n$, we have
\[
	1 - \tfrac{\ell}{n-j}
=
	1 - \tfrac \ell n \tfrac1{1-j/n}
=
	1 - \tfrac \ell n + \tfrac \ell n \tfrac{j/n}{1-j/n}
=
\rbb{ 1 - \tfrac \ell n } \rbb{ 1 + \Oh{ j \ell / n^2 } }.
\label{eq:k-cycle:product-approx}
\nt
\]

We turn first to $P_0$.
First note that $\defac{n-r-1}{k} = \defac{n-r}{k} \cdot \rbr{1 - k/(n-r)}$.
We have
\[
	\frac{\defac{n-r}{k}}{\defac nk}
=
	\frac{\defac{n-k}{r}}{\defac nr}
=
	\prodd[r-1]{s=0}
	\frac{n-k-s}{n-s}
=
	\prodd[r-1]{s=0}
	\rbbb{ 1 - \frac{k}{n-s} }.
\label{eq:k-cycle:P0-init}
\nt
\]
Combining (\ref{eq:k-cycle:product-approx}, \ref{eq:k-cycle:P0-init}),
and using the fact that $0 \le s \le r \asymp 1$,
we obtain
\[
	\defac{n-r}{k} / \defac nk
=
	\prodt[r-1]{s=0}
	(1 - k/n)
	\rbb{ 1 - \Oh{k/n^2} }
=
	\rbr{ 1 - k/n }^r
	\rbb{ 1 + \Oh{k/n^2} }.
\]
Now,
\(
	P_0
=
	\tfrac{n-r-1}{n-k} \cdot \defac{n-r}{k} / \defac nk.
\)
So, applying \eqref{eq:k-cycle:product-approx} again,
we obtain
\[
	P_0
=
	\rbr{ 1 - k/n }^{r+1}
	\rbb{ 1 + \Oh{k^2/n^2} }.
\label{eq:k-cycle:P0-final}
\nt
\]

We now turn to $P_1$ and $P_2$.
Using the approximation to $1 - \ell / (n-j)$, ie \eqref{eq:k-cycle:product-approx},
the following~hold:
\[
	P_1
&=
	\prodt[m]{i=2}
	\rbr{ 1 - k/n }
\cdot
	\rbb{ 1 + \Oh{k/n^2} }
=
	\rbr{ 1 - k/n }^{m-1} \rbb{ 1 + \Oh{k/n^2} };
\label{eq:k-cycle:P1-init}
\nt
\\
	P_2^{-1}
&=
	\prodt[m]{i=1}
	\rbr{ 1 - k/n }
\cdot
	\rbb{ 1 + \Oh{k/n^2} }
=
	\rbr{ 1 - k/n }^m \rbb{ 1 + \Oh{k/n^2} }.
\label{eq:k-cycle:P2-init}
\nt
\]
This uses $\max\bra{\lambda'_1, \bra{\lambda_i,\lambda'_i}_2^m, m} \le r \asymp 1$.
(Recall that $\lambda = (a_1, ..., a_m \mid b_1, ..., b_m)$.)
Hence
\[
	P_1 P_2
=
	\rbr{ 1 - k/n }^{-1} \rbb{ 1 + \Oh{k/n^2} }.
\label{eq:k-cycle:P12-final}
\nt
\]

Combining the expressions for $P_0$ and $P_1 P_2$,
ie (\ref{eq:k-cycle:P0-final}, \ref{eq:k-cycle:P12-final}),
we obtain
\[
	P_0 P_1 P_2
=
	\rbr{ 1 - k/n }^r \rbb{ 1 + \Oh{k/n^2} }.
\]
This is the main contribution to $s_\lambda(k)$; it remains to control the error in \cref{res:k-cycle:5a}.

If $k \gg 1$, then we necessarily have $r < k$, and so the error term is 0; if $k \asymp 1$, then the error term is $\Oh{n^{-k}} = \Oh{n^{-2}}$, as $k \ge 2$.
But $(1 - r/n)^k \asymp 1$, since $k \le \tfrac13 n$ and $r \asymp 1$, so this additive $\Oh{1/n^2}$ error is absorbed into the larger $\Oh{k/n^2}$ error.

In summary,
we have shown the desired expression for the character ratio $s_\lambda(k)$:
\[
	s_\lambda(k)
=
	\rbr{ 1 - r/n }^k
	\rbb{ 1 + \Oh{k/n^2} }.
\qedhere
\]
\end{Proof}

\begin{Proof}[Proof of \cref{res:k-cycle:short:k-big}]
	Choose $\theta' \cq 0.68$; then $\theta' - \tfrac16 > \tfrac12 + \tfrac1{10}$.
	Noting that
	\(
		n = \exp{\log n} \le \exp{\tfrac16k},
	\)
	inspection of the proof of \cite[Theorem~5(b)]{H:cutoff-k-cycle} gives the upper bound
	\[
		\abs{s_\lambda(k)}
	\le
		\expb{ k \rbb{ - \theta' + \tfrac16 + \oh1 } }
	\le
		\expb{ - (\tfrac12 + \tfrac1{10}) k }.
	\qedhere
	\]
\end{Proof}

\begin{Proof}[Proof of \cref{res:k-cycle:short:k-small}]
Under the given assumptions, \cite[Lemma~14]{H:cutoff-k-cycle} states that
\[
	\abs{s_\lambda(k)}
\le
	\rbb{ \sumt{a_i > k \sqrt n} a_i^k/n^k + \sumt{b_i > k \sqrt n} b_i^k/n^k }
	\rbb{ 1 + \Oh{\log n/n^{1/4}} }
+	\Ohb{ e^{-k} (\log n)^4 / n^{1/4} }.
\]
Further, we claim that if $r \cq n - \lambda_1$ satisfies $\tfrac13 n \le r \le n$ then
\[
	\sumt{} a_i^k/n^k + \sumt{} b_i^k/n^k
\le
	\expb{ - \tfrac{kr}{n} \rbr{\tfrac12 + \tfrac1{10} + \tfrac1{200}} }.
\]
Note that $1 + \Oh{\log n/n^{1/4}} \le \exp{\Oh{n^{-1/5}}}$ and $\Oh{e^{-k} (\log n)^4/n^{1/4}} = \oh{e^{-k}}$.
Since $r \asymp n$, and so $k r/n \asymp k \gtrsim 1$, these error terms can be absorbed by $\exp{ - \tfrac1{200} kr/n }$.
\cref{res:k-cycle:short:k-small} then follows.

It remains to prove our claim,
which is a slight sharpening of \cite[Lemma~15]{H:cutoff-k-cycle}.
The following claim comes from inspecting the proof of \cite[Lemma~15]{H:cutoff-k-cycle}:
	in order to prove that
	\[
		\sumt{} a_i^k/n^k + \sumt{} b_i^k/n^k
	\le
		\exp{ - ckr/n }
	\Quad{where}
		c \in \mbr,
	\]
	it suffices, writing $\delta \cq r/n \in [\tfrac13, 1]$, to prove that
	\[
		(1 - \delta)^{k-1} \le e^{-kc \delta},
	\Quad{ie}
		1 - \delta \le e^{-c \delta k/(k-1) }.
	\]
The worst case is clearly $k = 2$, in which case $k/(k-1) = 2$.
Thus we need $1 - \delta \le e^{-2c\delta}$.
If one can allow $\delta$ all the way down to 0, then one must take $c \le \tfrac12$; however, we only need $\delta \in [\tfrac13, 1]$.
One can then check that it is then sufficient to take $c$ so that
\[
	\tfrac23 = 1 - \tfrac13 \le e^{-2c/3},
\Quad{ie}
	c \le \tfrac32 \log(\tfrac32) \approx 0.608.
\]
In particular, we may take $c \cq \tfrac12 + \tfrac1{10} + \tfrac1{200} = 0.605$.
\end{Proof}

\section{Random Walks on Homogenous Spaces}
\label{sec:hom}

Throughout this section, $G$ will be a finite group and $K$ a subgroup.
Denote the \textit{homogenous space} consisting of the (right) cosets by
\(
	X \cq G/K \cq \bra{ g K \mid g \in G }.
\)
Denote the set of complex-valued functions on $X$ by
\(
	L(X) \cq \bra{ f : X \to \mbc }.
\)
We frequently identify this with the space of $K$ invariant functions on $G$, ie those $f : G \to \mbc$ for which $f(gk) = f(g)$ for all $g \in G$ and all $k \in K$.

\subsection{Gelfand Pairs and Spherical Fourier Analysis for Invariant Random Walks}
\label{sec:hom:gen}

The majority of this subsection---namely, the analysis leading up to \cref{res:hom:tv:exact}---is an abbreviated exposition of \cite[\S4]{CsST:harmonic-analysis-finite-groups}; a related exposition can be found in \cite[\S2]{CsST:gelfand-applications}.

Let $G$ be a finite group and let $K$ be a subgroup.
A function $f : G \to \mbc$ is $K$ \textit{bi-invariant} if
\[
	f(k_1 g k_2) = f(g)
\Quad{for all}
	g \in G
\Quad{and}
	k_1, k_2 \in K.
\]

\begin{defn}
\label{def:hom:gelfand}
	Let $G$ be a finite group and $K$ be a subgroup.
	The pair $(G, K)$ is called a \textit{Gelfand pair} if the algebra of $K$ bi-invariant functions (under convolution) is commutative.
	
	Equivalently, $(G, K)$ is a Gelfand pair if the permutation representation $\lambda$ of $G$ on $X$ defined by
	\(
		\rbr{ \lambda(g) f }(x) \cq f(g^{-1} x)
	\)
	for $g \in G$, $f \in L(X)$ and $x \in X$,
	is multiplicity-free.
\end{defn}

This equivalence is shown in \cite[Theorem~4.4.2]{CsST:harmonic-analysis-finite-groups}.
From now on, assume that $(G, K)$ is a Gelfand pair.
We next introduce \textit{spherical functions} and \textit{spherical representations}.

\begin{defn}
\label{def:hom:spherical:fn-rep}
	A $K$ bi-invariant function $\varphi : G \to \mbc$ is said to be \textit{spherical} if $\varphi(\id_G) = 1$ and
	\(
		\varphi * f
	=
		\rbb{ (\varphi * f)(\id_G) } \varphi
	\Quad{for all}
		\text{$K$ bi-invariant functions $f : G \mapsto \mbc$}.
	\)
	For a spherical function $\varphi$,
	the subspace of $L(X)$ generated by the $G$-translates of $\varphi$, ie
	\(
		V_\varphi
	\cq
		\langle \lambda(g) \varphi \mid g \in G \rangle
	\)
	where $\lambda$ is the \textit{permutation representation} of $G$ on $X$,
	is called the \textit{spherical representation}.
\end{defn}

For a representation $(\rho, V)$, write
\(
	V^K \cq \bra{ v \in V \mid \rho(k) v = v \: \forall \, k \in K }
\)
for the space~of~$K$~\textit{invariant} \textit{vectors} in $V$.
The following theorem is a culmination of statements from~\cite[\S4.5~and~\S4.6]{CsST:harmonic-analysis-finite-groups}.

\begin{thm}
\label{res:hom:spherical:res}
	The number of distinct spherical functions equals the number of orbits of $K$ on $X$.
	Denote by $\bra{\varphi_i}_0^N$ the distinct spherical functions, with $\varphi_0$ the constant function 1.

	Then
		\(
			L(X) = \oplus_0^N V_{\varphi_i},
		\)
		which is a multiplicity-free decomposition into irreps.
	Further, $\bra{\varphi_i}_0^N$ forms an orthogonal basis for the set of $K$ bi-invariant functions on $G$ with normalisation given, for each $i$, by
		\(
			\sumt{x \in X} \abs{\varphi_i(x)}^2
		=
			\abs{X} / d_i
		\)
		where $d_i \cq \dim V_{\varphi_i}$ is the dimension of the irrep $V_{\varphi_i}$.
	
	For any irrep $V$ we have $\dim V^K \le 1$ and $\dim V^K = 1$ if and only if $V$ is spherical.
\end{thm}

This allows us to construct a `spherical basis' in which the Fourier transform has a simple form.

\begin{defn}
\label{def:hom:spherical:ft}
	The \textit{spherical Fourier transform} $\widetilde \mu$ of a $K$ invariant function $\mu \in L(X)$ is defined~by%
	\[
		\widetilde \mu(i)
	\cq
		\sumt{x \in X}
		\mu(x) \overline{\varphi_i(x)}
	\Qfor
		i \in \bra{0, 1, ..., n}.
	\]
\end{defn}

\begin{cor}
	There exists an orthonormal basis of $K$ invariant functions on $G$ with the following property.
	Let $\mu$ be a $K$ bi-invariant function on $G$.
	If $(\tau, W)$ is a non-spherical irrep, then $\hat \mu(\tau) = 0$.
	If $(\rho_i, V_{\varphi_i})$ is a spherical irrep (with $i \in \bra{0, ..., N}$), then the matrix representing the operator $\widehat \mu(\rho_i)$ has only one non-zero entry, which is in the first position and has value $\abs K \widetilde \mu(i)$.
	
	As a consequence, a Fourier inversion formula holds:
	\[
		\mu^{*t}(x)
	=
		\abs{X}^{-1}
		\sumt[N]{i=0}
		d_i \varphi_i(x) \widetilde \mu(i)^t
	\Qforall
		x \in X
	\Qand
		t \in \mbn_0,
	\]
	where $\mu^{*t}$ is the $t$-fold self-convolution of $\mu$.
\end{cor}

From this we immediately obtain for the TV distance between $\mu^{*t}$ and $\Unif_X$.
To apply this to random walks on $G$, the step distribution must be $K$ bi-invariant;
this is the case if the stochastic transition matrix $P = (p_{x,y})_{x,y \in X}$ is \textit{$G$-invariant}:
\(
	p_{x,y} = p_{gx,gy}
\)
for all $x,y \in X$ and all $g \in G$.

When looking at such random walks, we always start from a point which is stabilised by $K$.

\begin{defn}
	Let $G$ be a finite group and $K$ be a subgroup.
	Let $G$ act on the homogenous space $X \cq G/K$ by the left coset action:
	\(
		g \cdot (hK) \cq (gh) K.
	\)
	Say $\bar x \in K$ is \textit{stabilised by $K$} if $k \cdot \bar x = \bar x$ for all $k \in K$. Equivalently, $\bar g K$ is \textit{stabilised by $K$} if and only if $K = \bar g K \bar g^{-1}$.
\end{defn}

When starting a random walk with $G$-invariant transition matrix from $\bar x \in X = G/K$ which is stabilised by $K$,
one can then check $P^t(\bar x, \cdot) = \mu_{\bar x}^{*t}(\cdot)$ for all $t \in \mbn_0$ where $\mu_{\bar x}(\cdot) \cq P(\bar x, \cdot)$;
that is, the probability of being at $x$ after $t$ steps when started from $\bar x$ is $\mu_{\bar x}^{*t}(x)$ for all $x \in X$ and all $t \in \mbn_0$.
Altogether, we have now proved the following proposition.

\begin{prop}[{\cite[Proposition~4.9.1]{CsST:harmonic-analysis-finite-groups}}]
\label{res:hom:tv:exact}
	Let $(G, K)$ be a Gelfand pair and denote $X \cq G/K$.
	Let $\bra{\varphi_i}_{i=0}^N$ be the associated spherical functions, considered as $K$ bi-invariant functions on $X$, and $\bra{d_i}_{i=0}^N$ the associated dimensions; assume that $\varphi_0(x) = 1$ for all $x \in X$.
	
	Let $\bar x$ be an element of $X$ stabilised by $K$.
	Let $P$ be a $G$-invariant stochastic matrix and set $\mu_{\bar x}(\cdot) \cq p_{\bar x, \cdot}$.
	Let $t \in \mbn_0$ and $x \in X$.
	Then
	\[
		\mu_{\bar x}^{*t}(x) - \abs X^{-1}
	=
		\abs X^{-1}
		\sumt[N]{i=1}
		d_i \varphi_i(x) \widetilde \mu_{\bar x}(i)^t,
	\]
	where $\widetilde \mu_{\bar x}$ is the spherical Fourier transform of $\mu_{\bar x}$.
	As a corollary,
	we have
	\[
		d_\TV\rbb{ P^t(\bar x, \cdot), \: \Unif_X }
	=
		\tfrac12 \abs X^{-1} \sumt{x \in X} \absb{ \sumt[N]{i=1} d_i \varphi_i(x) \widetilde \mu_{\bar x}(i)^t }.
	\]
\end{prop}

We now have all the ingredients to prove our TV-approximation lemma for random walks on homogeneous spaces corresponding to Gelfand pairs, ie \cref{res:tv-approx:hom}; we rested it here for convenience.

\begin{lem}[TV Approximation Lemma]
\label{res:hom:tv:approx}
	Let $(G, K)$ be a Gelfand pair and denote $X \cq G/K$.
	Let $\bar x$ be an element of $X$ stabilised by $K$.
	Let $\bra{\varphi_i}_{i=0}^N$ be the associated spherical functions, considered as $K$ bi-invariant functions on $X$, and $\bra{d_i}_{i=0}^N$ the associated dimensions; assume that $\varphi_0(x) = 1$ for all $x \in X$.
	Let $P$ be a $G$-invariant stochastic matrix and set $\mu_{\bar x}(\cdot) \cq P\rbr{\bar x, \cdot}$.
	
	Let $t \in \mbn_0$ and $I \subseteq \bra{1, ..., N}$.
	Then
	\[
		\absB{
			d_\TV\rbb{ P^t(\bar x, \cdot), \: \Unif_X }
		-	\tfrac12 \abs X^{-1} \sumt{x \in X} \absb{ \sumt{i \in I} d_i \varphi_i(x) \widetilde \mu_{\bar x}(i)^t }
		}
	\le
		\tfrac12 \sumt{i \notin I} \sqrt{d_i} \abs{ \widetilde \mu_{\bar x}(i) }^t,
	\]
	where $\widetilde \mu_{\bar x} : i \mapsto \sumt{x \in X} \mu_{\bar x}(x) \overline{\varphi_i(x)}$ is the spherical Fourier transform of $\mu_{\bar x}$.
\end{lem}

\begin{Proof}
First we apply \cref{res:hom:tv:exact} and the triangle inequality:
\[
&	\absB{
		d_\TV\rbr{ \mu_{\bar x}^{*k}, \: \pi }
	-	\tfrac12 \abs X^{-1} \sumt{x \in X} \absb{ \sumt{i \in I} d_i \varphi_i(x) \widetilde \mu_{\bar x}(i)^k }
	}
\le
	\tfrac12 \abs X^{-1} \sumt{x \in X} \absb{ \sumt{i \notin I} d_i \varphi_i(x) \widetilde \mu_{\bar x}(i)^k }
\\&\qquad
\le
	\tfrac12 \abs X^{-1} \sumt{x \in X} \sumt{i \notin I} d_i \abs{ \varphi_i(x) \widetilde \mu_{\bar x}(i)^k }
=
	\tfrac12 \sumt{i \notin I} d_i \abs{\widetilde \mu_{\bar x}(i)}^k \cdot \abs X^{-1} \sumt{x \in X} \abs{\varphi_i(x)}.
\]
Applying Cauchy--Schwarz and the standard spherical orthogonality relations (see, eg, \cite[Proposition~4.7.1]{CsST:harmonic-analysis-finite-groups} or \cite[Equation~(2.11)]{CsST:gelfand-applications}), we obtain
\[
	\rbb{ \sumt{x \in X} \abs{\varphi_i(x)} }^2
\le
	\abs X \sumt{x \in X} \abs{\varphi_i(x)}^2
=
	\abs X \cdot \abs X / d_i.
\]
Plugging this into the previous bound,
we deduce the lemma.
\end{Proof}

\subsection{Limit Profile for Many-Urn Ehrenfest Diffusion}
\label{sec:hom:e-urn}

Suppose that one has $\nn$ \textit{balls} labelled $1$ through $\nn$ and $m+1$ \textit{urns} labelled $0$ through $m$.
The set of all configurations can be identified with the set
\(
	X_{\nn,m+1} \cq \bra{0, 1, ..., m}^\nn:
\)
an element $x = (x_1, ..., x_\nn) \in X_{\nn,m+1}$ indicates that the $j$-th ball is in the $x_j$-th urn.
Initially, put all the balls in the first urn (labelled $0$):
	this is the \textit{initial configuration}, and corresponds to $\bar x \cq (0, 0, ..., 0)$.

We can endow $X$ with a metric structure:
	for $x,y \in X_{\nn,m+1}$,
	set
	\[
		d(x,y) \cq \absb{ \bra{ k \in [\nn] \mid x_k \ne y_k } }.
	\]
Thinking of $x$ and $y$ as configurations of balls, $d(x,y)$ is the number of balls which are not in the same urn in the two configurations.

We consider the random walk on $X \cq X_{\nn,m+1}$ described by the following step:
	choose uniformly at random a ball and an urn;
	put the chosen ball in the chosen urn.
In terms of a transition matrix $R$ on $X \times X$, this is given by the following expressions, for $x,y \in X$:
\[
	R(x,y) = \tfrac1{m+1} \text{ if } x = y;
\quad
	R(x,y) = \tfrac1{\nn(m+1)} \text{ if } d(x,y) = 1;
\quad
	R(x,y) = 0 \ \text{otherwise}.
\]

The following theorem is a restatement of \cref{res:intro:e-urn}, but written more formally: cutoff is for a sequence of Markov chains; we make this sequence explicit.

\renewcommand{\cc}{\ensuremath{c}}
\begin{thm}[Limit Profile for Generalised Ehrenfest Urn]
\label{res:hom:e-urn:res}
	Let $\nn,m \in \mbn$.
	Consider $\nn$ balls labelled $1, ..., \nn$ and $m+1$ urns labelled $0, 1, ..., m$.
	Consider the following Markov chain:
		at each step,
		choose a ball and an urn uniformly and independently;
		place said ball in said urn.
	For $\tt \in \mbn_0$, write $d_\TV^{\nn,m}(t)$ for the TV distance of this Markov chain after $t$ steps from its invariant distribution when started with all $n$ balls initially in the urn labelled 0.
	
	Let $(\nn_N)_\Ninn, (\mm_N)_\Ninn \in \mbn^\mbn$.
	Suppose that $\lim_N \mm_N / \nn_N = 0$.
	Then,
	for all $\cc \in \mbr$,
	we have
	\[
		d_\TV^{\nn_N,\mm_N}\rbb{ \tfrac12 \nn_N \log(\mm_N \nn_N) + \cc \nn_N }
	\to
		2 \, \Phi\rbb{ \tfrac12 e^{-\cc} } - 1
	\quad
		\asinf N.
	\]
\end{thm}

As in previous sections, for ease of presentation we omit the $N$-subscripts in the proof.
We start by phrasing the Ehrenfest urn model in Gelfand pair language.
To do this, we give a very abbreviated exposition of \cite[\S3]{CsST:gelfand-applications}.
Let $\symgr_{m+1}$ and $\symgr_n$ be the symmetric groups on $\bra{0, 1, ..., m}$ and $\bra{1, ..., n}$, respectively.
Then $X_{n,m+1} = \bra{0, 1, ..., n}^n$ is a homogenous space for the wreath product $\symgr_{m+1} \wr s_n$ under the action
\(
	\rbr{ \sigma_1, ..., \sigma_n; \theta } \cdot \rbr{ x_1, ..., x_n }
\cq
	\rbr{ \sigma_1 x_{\theta^{-1}(1)}, ..., \sigma_n x_{\theta^{-1}(n)} },
\)
ie $X_i$ is moved by $\theta$ to the position $\theta(i)$ and then it is changed by the action of $\sigma_{\theta(i)}$.
Note that the stabiliser of $\bar x \cq (0, 0, ..., 0) \in X_{n,m+1}$ coincides with the wreath product $\symgr_m \wr \symgr_n$, where $\symgr_m \le \symgr_{m+1}$ is the stabiliser of $0$.
Therefore we can write
\(
	X_{m+1,n}
=
	\rbr{ \symgr_{m+1} \wr \symgr_n } / \rbr{ \symgr_m \wr \symgr_n }.
\)
The action is distance transitive, and so the group $\symgr_{m+1} \wr \symgr_n$ acts isometrically on $X_{n,m+1}$.
It follows that $\rbr{ \symgr_{m+1} \wr \symgr_n } / \rbr{ \symgr_m \wr \symgr_n }$ is a Gelfand pair; see \cite[Example~2.5]{CsST:gelfand-applications}.

The associated spherical functions and dimensions are given by the following proposition.

\begin{thm}[Spherical Functions; {\cite[Theorem~3.1]{CsST:gelfand-applications}}]
\label{res:hom:e-urn:spherical-fns}
	For each $i \in \bra{0, 1, ..., n}$,
		the dimension $d_i$ satisfies
		\(
			d_i
		=
			m^i \binomt ni
		\)
	and
		the spherical function $\varphi_i$ satisfies
		\[
			\varphi_i(x)
		=
			\binom ni^{-1}
			\sumd[\min\bra{\ell,i}]{r=\max\bra{0,i-n-\ell}}
			\binom \ell r \binom{n-\ell}{i-r} \rbbb{-\frac1m}^r
		\Qfor
			x \in X
		\Qwhere
			\ell \cq d(\bar x,x)
		\]
\end{thm}

\begin{rmkt}
	The spherical functions are the \textit{Krawtchouk polynomials}, given in \cref{def:rev:gibbs:adj}:
	\[
		\varphi_i(x)
	\cq
		\varphi_i(\ell)
	\cq
		K_i\rbb{ \ell; \, \tfrac{m}{m+1}, \, n }
	\Quad{where}
		\ell \cq d(\bar x,x),
	\]
	using the notation there.
	These are orthogonal with respect to the Binomial measure by \cref{res:rev:gibbs:kraw-orthog}:
	\[
		\sumt[n]{\ell=0}
			K_i\rbb{ \ell; \, \tfrac{m}{m+1}, \, n} K_j\rbb{ \ell; \, \tfrac{m}{m+1}, \, n}
		\cdot
			\binomt n\ell m^\ell/(m+1)^n
	=
		\rbb{ m^i \binomt ni }^{-1} \delta_{i,j}.
	\]
	This can also be seen as consequence of the orthogonality of spherical functions, ie \cref{res:hom:spherical:res}.
\end{rmkt}


We first determine the spherical Fourier transform of the step distribution $\mu(\cdot) \cq R(\bar x, \cdot)$.

\begin{lem}[Spherical Fourier Transform]
\label{res:hom:e-urn:spherical-fourier}
	For all $\ii \in \bra{0, 1, ..., \nn}$,
	we have
	\(
		\widetilde \mu(\ii)
	=
		1 - \ii/\nn.
	\)
\end{lem}

\begin{Proof}
Noting the slight laziness,
we have
\[
	\widetilde \mu(\ii)
=
	\tfrac{m}{m+1} \rbb{ \tfrac1m + \varphi_\ii(1) }.
\]
Using the expression for $\varphi_\ii(1)$ given by \cref{res:hom:e-urn:spherical-fns},
we obtain
\(
	\widetilde \mu(\ii)
=
	1 - \ii/n.
\)
	%
\end{Proof}

There are $m^\ell \binom \nn\ell$ different $x$ with $d(\bar x,x) = \ell$.
Applying \cref{res:hom:e-urn:spherical-fns},
we obtain the following expressions for the terms in \cref{res:tv-approx:hom}:
\[
	\MT
&
\cq
	\abs X^{-1} \sumt{x \in X}
	\absb{ \sumt[\MM]{\ii=1} d_\ii \varphi_\ii(x) \widetilde \mu(\ii)^\tt }
\\&
=
	(m+1)^{-\nn} \sumt[\nn]{\ell=0} m^\ell \binomt \nn\ell
	\absb{ \sumt[\MM]{\ii=1} m^\ii \binomt \nn\ii \rbr{ 1 - \ii/\nn }^\tt \varphi_\ii(\ell) };
\\
	\ET
&
\cq
	\sumt{\ii \notin \II} \sqrt{d_\ii} \abs{ \widetilde \mu(\ii) }^\tt
=
	\sumt{\ii > \MM}
	\mm^{\ii/2} \binomt \nn\ii^{1/2} \rbr{1 - \ii/\nn}^\tt.
\]

Our first aim is to use this to determine which are the `important' spherical statistics.

\begin{lem}[Error Term]
\label{res:hom:e-urn:error}
	For all $\eps > 0$ and all $\cc \in \mbr$,
	there exists an $\MM \cq \MM(\cc,\eps)$ so that,
	for $\tt \cq \tfrac12 n \log(mn) + \cc n$,
	if $\II \cq \bra{1, ..., M}$, then
	\[
		\ET
	\le
		\ET'
	\le
		\eps
	\Qwhere
		\ET'
	\cq
		\sumt{\ii \notin \II}
		\sqrt{d_\ii} e^{-\ii \tt/\nn}
	=
		\sumt{\ii > \MM}
		\mm^{\ii/2} \binomt \nn\ii^{1/2} e^{-\ii \tt/\nn}
	\le
		\eps.
	\]
\end{lem}

\begin{Proof}
Using \cref{res:hom:e-urn:spherical-fourier}, we have
\(
	\abs{\widetilde \mu(\ii)} \le e^{-\ii/\nn}
\)
for all $\ii$.
The inequality $\ET \le \ET'$ now follows.
The equality in the definition of $\ET'$ is now an immediate consequence of \cref{res:hom:e-urn:spherical-fns}.
For the inequality $\ET' \le \eps$,
choose $\MM$ so that
\(
	\sumt{\ii > \MM}
	e^{-\cc\ii} / \sqrt{\ii!}
\le
	\eps.
\)
Then we have
\[
	\ET
\le
	\sumt{\ii > \MM}
	\rbb{ (\mm \nn)^{\ii/2} e^{-\tt/\nn} }^\ii / \sqrt{\ii!}
=
	\sumt{\ii > \MM}
	e^{-\cc\ii} / \sqrt{\ii!}
\le
	\eps.
\qedhere
\]
\end{Proof}

From now on, choose $\MM \cq \MM(\cc,\eps)$ as in \cref{res:hom:e-urn:error}.
Hence, for the main term, we need only deal with spherical statistics with $i \asymp 1$.
We would then like to use the replacement $\lambda_\ii \approx e^{-\ii/\nn}$.

\begin{defn}[Adjusted Main Term]
\label{def:hom:e-urn:main:adj}
	Recalling that $\tt = \tfrac12 \nn \log(\mm \nn) + \cc \nn$,
	define
	\[
		\MT'
	\cq
		(m+1)^{-\nn} \sumt[\nn]{\ell=0} m^\ell \binomt \nn\ell
		\absb{ \sumt{i\ge1} \binomt \nn\ii \varphi_i(\ell) \mm^{\ii/2} e^{-\cc\ii} / \nn^{\ii/2} }.
	\]
\end{defn}

Conveniently, the adjusted main term $\MT'$ in this case (\cref{def:hom:e-urn:main:adj}) is \textit{exactly} the same as that for the Gibbs sampler (see \cref{def:rev:gibbs:adj}) in \S\ref{sec:rev:gibbs}; to match notation, replace $\mm$ with $\aa$.

The following two lemmas are simply a restatement of \cref{res:rev:gibbs:main:approx,res:rev:gibbs:main:eval}.

\begin{subtheorem}{thm}
	\label{res:hom:e-urn:main}

\begin{lem}[Main Term: Approximation]
\label{res:hom:e-urn:main:approx}
	For all $\eps > 0$ and all $\cc \in \mbr$,
	with $\MM \cq \MM(\cc,\eps)$,
	we~have
	\[
		\absb{ \MT - \MT' }
	\le
		2 \eps.
	\]
\end{lem}

\begin{lem}[Main Term: Evaluation]
\label{res:hom:e-urn:main:eval}
	For all $\cc \in \mbr$,
	with $\MM \cq \MM(\cc,\eps)$,
	we have
	\[
		\tfrac12 \MT'
	\to
		2 \, \Phi\rbb{ \tfrac12 e^{-\cc} } - 1.
	\]
\end{lem}

\end{subtheorem}

We now have all the ingredients to establish the limit profile for the Ehrenfest urn model.

\begin{Proof}[Proof of \cref{res:hom:e-urn:res}]
Let us summarise what we have proved.
These are all evaluated at the target mixing time
\(
	\tt
=
	\tfrac12 \log(\mm \nn) + \cc \nn
\)
with $\MM \cq \MM(\cc,\eps)$ given by \cref{res:hom:e-urn:error}.
\begin{itemize}[noitemsep, topsep = \smallskipamount, label = \bcdot]
	\item 
	By \cref{res:hom:e-urn:error},
	the error term $\ET$ satisfies
	\(
		\ET \le \eps.
	\)
	
	\item 
	By \cref{res:hom:e-urn:main:approx},
	the original main term $\MT$ satisfies
	\(
		\abs{\MT - \MT'} \le 2 \eps.
	\)
	
	\item 
	By \cref{res:hom:e-urn:main:eval},
	the adjusted main term $\MT'$ satisfies
	\(
		\tfrac12 \MT' \to 2 \, \Phi(\tfrac12 e^{-\cc}) - 1
	\)
	\asinf \nn.
\end{itemize}
Since $\eps > 0$ is arbitrary,
applying the TV-approximation lemma for random walks on homogenous spaces, namely \cref{res:tv-approx:hom},
we immediately deduce the theorem.
\end{Proof}

\section*{References}

\renewcommand{\bibfont}{\sffamily}
\renewcommand{\bibfont}{\sffamily\small}
\printbibliography[heading=none]

@incollection{A:rw-group-mixing,
  title = {Random Walks on Finite Groups and Rapidly Mixing Markov Chains},
  booktitle = {Seminar on Probability, XVII},
  author = {Aldous, David},
  date = {1983},
  volume = {986},
  pages = {243–297},
  publisher = {Springer, Berlin},
  doi = {10.1007/BFb0068322},
  url = {https://doi.org/10.1007/BFb0068322},
  mrclass = {60J15 (60B15)},
  mrnumber = {770418},
  mrreviewer = {A. Mukherjea},
  series = {Lecture Notes in Math.}
}

@article{BhS:cutoff-nonbacktracking,
  title = {Cutoff for Nonbacktracking Random Walks on Sparse Random Graphs},
  author = {Ben-Hamou, Anna and Salez, Justin},
  date = {2017},
  journaltitle = {Ann. Probab.},
  volume = {45},
  pages = {1752-1770},
  issn = {0091-1798},
  doi = {10.1214/16-AOP1100},
  url = {https://doi.org/10.1214/16-AOP1100},
  fjournal = {The Annals of Probability},
  issue = {3},
  mrclass = {60J10 (05C80 05C81 60G50)},
  mrnumber = {3650414},
  mrreviewer = {Malwina Joanna Luczak},
  number = {3}
}

@article{BS:cutoff-conj-inv,
  title = {Cutoff for Conjugacy-Invariant Random Walks on the Permutation Group},
  author = {Berestycki, Nathanaël and Şengül, Batı},
  date = {2019},
  journaltitle = {Probab. Theory Related Fields},
  volume = {173},
  pages = {1197-1241},
  issn = {0178-8051},
  doi = {10.1007/s00440-018-0844-y},
  url = {https://doi.org/10.1007/s00440-018-0844-y},
  fjournal = {Probability Theory and Related Fields},
  issue = {3-4},
  mrclass = {60B15 (60C05 60G50)},
  mrnumber = {3936154},
  number = {3-4}
}

@article{BSZ:k-cycle,
  title = {Mixing Times for Random $k$-Cycles and Coalescence-Fragmentation Chains},
  author = {Berestycki, Nathanaël and Schramm, Oded and Zeitouni, Ofer},
  date = {2011},
  journaltitle = {Ann. Probab.},
  volume = {39},
  pages = {1815-1843},
  issn = {0091-1798},
  doi = {10.1214/10-AOP634},
  url = {https://doi.org/10.1214/10-AOP634},
  fjournal = {The Annals of Probability},
  issue = {5},
  mrclass = {60B15 (60J27)},
  mrnumber = {2884874},
  mrreviewer = {Sara Brofferio},
  number = {5}
}

@article{CsST:gelfand-applications,
  title = {Finite Gel’fand Pairs and Their Applications to Probability and Statistics},
  author = {Ceccherini-Silberstein, Tullio and Scarabotti, Fabio and Tolli, Filippo},
  date = {2007-02-01},
  journaltitle = {Journal of Mathematical Sciences},
  shortjournal = {Journal of Mathematical Sciences},
  volume = {141},
  pages = {1182-1229},
  issn = {1573-8795},
  doi = {10.1007/s10958-007-0041-5},
  url = {https://doi.org/10.1007/s10958-007-0041-5},
  issue = {2},
  number = {2}
}

@book{CsST:harmonic-analysis-finite-groups,
  title = {Harmonic Analysis on Finite Groups},
  author = {Ceccherini-Silberstein, Tullio and Scarabotti, Fabio and Tolli, Filippo},
  date = {2008},
  volume = {108},
  publisher = {Cambridge University Press, Cambridge},
  doi = {10.1017/CBO9780511619823},
  url = {https://doi.org/10.1017/CBO9780511619823},
  isbn = {978-0-521-88336-8},
  mrclass = {43-01 (20C05 20C30 33D80 43A90 60J10 60J20)},
  mrnumber = {2389056},
  mrreviewer = {Alain Valette},
  pagetotal = {xiv+440},
  series = {Cambridge Studies in Advanced Mathematics}
}

@book{D:group-rep,
  title = {Group Representations in Probability and Statistics},
  author = {Diaconis, Persi},
  date = {1988},
  volume = {11},
  publisher = {Institute of Mathematical Statistics, Hayward, CA},
  isbn = {0-940600-14-5},
  mrclass = {60-02 (20C99 62-02)},
  mrnumber = {964069},
  mrreviewer = {Philippe Bougerol},
  pagetotal = {vi+198},
  series = {Institute of Mathematical Statistics Lecture Notes---Monograph Series}
}

@article{DGM:rw-hypercube,
  title = {Asymptotic Analysis of a Random Walk on a Hypercube with Many Dimensions},
  author = {Diaconis, Persi and Graham, R. L. and Morrison, J. A.},
  date = {1990},
  journaltitle = {Random Structures Algorithms},
  volume = {1},
  pages = {51-72},
  issn = {1042-9832},
  doi = {10.1002/rsa.3240010105},
  url = {https://doi.org/10.1002/rsa.3240010105},
  fjournal = {Random Structures & Algorithms},
  issue = {1},
  mrclass = {60J10 (60C05)},
  mrnumber = {1068491},
  mrreviewer = {David J. Aldous},
  number = {1}
}

@article{DKSc:exp-families,
  title = {Gibbs Sampling, Exponential Families and Orthogonal Polynomials},
  author = {Diaconis, Persi and Khare, Kshitij and Saloff-Coste, Laurent},
  date = {2008-05},
  journaltitle = {Statistical Science. A Review Journal of the Institute of Mathematical Statistics},
  shortjournal = {Statist. Sci.},
  volume = {23},
  pages = {151-178},
  issn = {0883-4237},
  doi = {10.1214/07-STS252},
  url = {https://doi.org/10.1214/07-STS252},
  fjournal = {Statistical Science. A Review Journal of the Institute of Mathematical Statistics},
  mrclass = {60E05 (60J10 62D05)},
  mrnumber = {2446500},
  mrreviewer = {Mu Fa Chen},
  number = {2}
}

@article{DS:bernoulli-laplace,
  title = {Time to Reach Stationarity in the Bernoulli--Laplace Diffusion Model},
  author = {Diaconis, Persi and Shahshahani, Mehrdad},
  date = {1987},
  journaltitle = {SIAM J. Math. Anal.},
  volume = {18},
  pages = {208-218},
  issn = {0036-1410},
  doi = {10.1137/0518016},
  url = {https://doi.org/10.1137/0518016},
  fjournal = {SIAM Journal on Mathematical Analysis},
  issue = {1},
  mrclass = {60B15 (60J10)},
  mrnumber = {871832},
  mrreviewer = {Philippe Bougerol},
  number = {1}
}

@article{DS:random-trans,
  title = {Generating a Random Permutation with Random Transpositions},
  author = {Diaconis, Persi and Shahshahani, Mehrdad},
  date = {1981},
  journaltitle = {Z. Wahrsch. Verw. Gebiete},
  volume = {57},
  pages = {159-179},
  issn = {0044-3719},
  doi = {10.1007/BF00535487},
  url = {https://doi.org/10.1007/BF00535487},
  fjournal = {Zeitschrift für Wahrscheinlichkeitstheorie und Verwandte Gebiete},
  issue = {2},
  mrclass = {60C05 (60J15)},
  mrnumber = {626813},
  mrreviewer = {Lars Holst},
  number = {2}
}

@article{EE:e-urn,
  title = {Über Zwei Bekannte Einwäande Gegen das Boltzmann\-sche H-Theorem},
  author = {Ehrenfest, Tatiana and Ehrenfest, Paul},
  date = {1907},
  journaltitle = {Physikalische Zeitschrift 8},
  pages = {311-314}
}

@article{H:cutoff-k-cycle,
  title = {The Random $k$ Cycle Walk on the Symmetric Group},
  author = {Hough, Robert},
  date = {2016},
  journaltitle = {Probab. Theory Related Fields},
  volume = {165},
  pages = {447-482},
  issn = {0178-8051},
  doi = {10.1007/s00440-015-0636-6},
  url = {https://doi.org/10.1007/s00440-015-0636-6},
  fjournal = {Probability Theory and Related Fields},
  issue = {1-2},
  mrclass = {60J10 (20C30 60B15)},
  mrnumber = {3500276},
  mrreviewer = {Sara Brofferio},
  number = {1-2}
}

@online{HOt:rcg:abe:cutoff,
  title = {Cutoff for Almost All Random Walks on Abelian Groups},
  author = {Hermon, Jonathan and Olesker-Taylor, Sam},
  date = {2021-02-04},
  url = {http://arxiv.org/abs/2102.02809},
  urldate = {2021-02-05},
  archiveprefix = {arXiv},
  eprint = {2102.02809},
  eprinttype = {arxiv},
  primaryclass = {math}
}

@article{KS:hypergeo-askey,
  title = {The Askey-Scheme of Hypergeometric Orthogonal Polynomials and Its $q$-Analogue},
  author = {Koekoek, Roelof and Swarttouw, René F.},
  date = {1998},
  journaltitle = {Delft University of Technology},
  online = {homepage.tudelft.nl/11r49/askey/}
}

@article{LP:ramanujan,
  title = {Cutoff on All Ramanujan Graphs},
  author = {Lubetzky, Eyal and Peres, Yuval},
  date = {2016},
  journaltitle = {Geom. Funct. Anal.},
  volume = {26},
  pages = {1190-1216},
  issn = {1016-443X},
  doi = {10.1007/s00039-016-0382-7},
  url = {https://doi.org/10.1007/s00039-016-0382-7},
  fjournal = {Geometric and Functional Analysis},
  issue = {4},
  mrclass = {60J10 (05C81)},
  mrnumber = {3558308},
  mrreviewer = {Serguei Popov},
  number = {4}
}

@book{LPW:markov-mixing,
  title = {Markov Chains and Mixing Times},
  author = {Levin, David A. and Peres, Yuval and Wilmer, Elizabeth L.},
  date = {2017},
  edition = {Second},
  publisher = {American Mathematical Society, Providence, RI, USA},
  doi = {10.1090/mbk/107},
  url = {https://doi.org/10.1090/mbk/107},
  isbn = {978-1-4704-2962-1},
  mrclass = {60J10 (60-01 60B15 60C05 60J27 60K35 68U20 82C22)},
  mrnumber = {3726904},
  pagetotal = {xvi+447},
  zmnumber = {1390.60001}
}

@article{S:mixing-notes,
  title = {Temps de Mélange des Chaînes de Markov (in French)},
  author = {Salez, Justin},
  date = {2018},
  journaltitle = {Online Lecture Notes},
  customeprint = {www.ceremade.dauphine.fr/~salez/mixing.pdf}
}

@article{T:limit-profile,
  title = {Limit Profile for Random Transpositions},
  author = {Teyssier, Lucas},
  date = {2020},
  journaltitle = {Annals of Probability},
  shortjournal = {Ann. Probab.},
  volume = {48},
  pages = {2323–2343},
  issn = {0091-1798},
  doi = {10.1214/20-AOP1424},
  url = {https://doi.org/10.1214/20-AOP1424},
  fjournal = {The Annals of Probability},
  mrclass = {60J10 (20C30 60B15)},
  mrnumber = {4152644},
  number = {5},
  zmnumber = {07276926}
}

\section{Appendix}

\subsection{Simple Random Walk on the Hypercube}
\label{sec:rw-hypercube}

Let $G = \mbf_2^n$ with identity $\id = (0, ..., 0) \in \mbf_2^n$ and $\mu$ be the probability measure, such that $\mu(\id) = 1/2$ and $\mu(e_i) = 1/(2n)$ for all $i \in [n]$, where $e_i$ is the vector in $G$ that has all entries equal to zero, but the $i$-th one which is equal to 1.
Then, for $t \in \mbn$, the law of the random walk with step distribution run for $t$ steps and started from $\id \in \mbf_2^n$ is given by the $t$-fold convolution $\mu^{*t}$.

Write $\Phi(\cdot)$ for the cdf of a standard normal distribution.
The following theorem gives the limit profile for the simple random walk on the hypercube.

\begin{thm}
\label{res:rw-hypercube:res}
	Let $n \in \mbn$ and consider the simple random walk on the hypercube $\mbf_2^n$.
	For $t \in \mbn_0$, write $d_\TV^n(t)$ for the TV distance of this Markov chain after $t$ steps from its invariant distribution.
	
	Then,
	for all $\cc \in \mbr$,
	we have
	\[
		d_\TV^n\rbb{ \tfrac12 n \log n + \cc n }
	\to
		2 \, \Phi\rbb{ \tfrac12 e^{-\cc} } - 1
	=
		\Phi\rbb{   \tfrac12 e^{-\cc} }
	-	\Phi\rbb{ - \tfrac12 e^{-\cc} }
	\quad
		\asinf n.
	\]
\end{thm}

\begin{Proof}
Observe that $d_\TV(\tt) = d_\TV(\mu^{*\tt}, \: \Unif_{\mbf_2^n})$ for all $\tt \in \mbf_2^n$ since $\mu$ is the step distribution of the simple random walk on the hypercube $\mbf_2^n$.

Since the group $\mbf_2^n$ is Abelian, the irreps are indexed by elements of the group.
Set
\[
	\chi_x(y)
\cq
	(-1)^{x \cdot y}
\Qwhere
	x \cdot y
\cq
	\sumt[n]{1} x_i y_i
\Qfor
	x,y \in \bra{0,1}^n,
\]
with addition modulo 2.
It is not difficult to check that these are the irreps (or equivalently characters as the group is Abelian so all irreps are of dimension 1); see, eg, \cite[\S2.3]{CsST:harmonic-analysis-finite-groups}.
Note that $(0, ..., 0) \in \mbf_2^n$ corresponds to the trivial partition.
Taking the Fourier transform, we see that
\[
	\widehat \mu(x)
=
	\sumt{y \in \mbf_2^n}
	(-1)^{x \cdot y} \mu(y)
=
	\tfrac12
+	\tfrac12 n^{-1} \rbb{ \abs{ \bra{ i \mid x_i = 0 } } - \abs{ \bra{ i \mid x_i = 1 } } }
=
	1 - \abs x / n
\]
where $\abs x \cq \sumt{i} x_i = \abs{ \bra{ i \mid x_i = 1 } }$ is the Hamming weight.
The Fourier inversion formula
gives
\[
	2 \, d_\TV\rbb{ \mu^{*\tt}, \: \Unif_{\mbf_2^n} }
=
	2^{-n}
	\sumt{x \in \mbf_2^n}
	\absb{
		\sumt{y \in \mbf_2^n \setminus \bra{(0,...,0)}}
		(-1)^{x \cdot y} \rbr{ 1 - \abs y / n }^\tt
	}
\Qforall
	\tt \in \mbn_0;
\]
see, eg, \cite[\S3.10]{CsST:harmonic-analysis-finite-groups}.
We now compute the inner sum.
Note that for each value of $\abs x$, there are $\binom n{\abs x}$ different $x \in \mbf_2^n$ which have this value.
By convention, set $\binom Nr \cq 0$ unless $0 \le r \le N$.~%
We~have
\[
	\sumt{ y \in \mbf_2^n \setminus \bra{(0,...,0)} }
	(-1)^{x \cdot y} \rbr{ 1 - \abs y/n }^\tt
&
=
	\sumt[n]{i=0}
	\rbr{ 1 - j/n }^\tt
	\sumt{i \ge 0}
	\binomt{\abs x}{i} \binomt{n - \abs x}{j - i}
	(-1)^i
-	1
\\&
=
	\sumt[\abs x]{i=0}
	\binomt{\abs x}{i}
	(-1)^i
	\sumt[n - \abs x + i]{j=i}
	\rbr{ 1 - j/n }^t
	\binomt{n - \abs x}{j - i}
-	1.
\]
Letting $\cc \in \mbr$ and setting
\(
	t \cq \tfrac12 n \log n + \cc n,
\)
we now have
\[
&	2 \cdot d_\TV\rbb{ \mu^{*t}, \: \Unif_{\mbf_2^n} }
=
	2^{-n}
	\sumt{x \in \mbf_2^n}
	\absb{
		\sumt[\abs x]{i=0}
		\binomt{\abs x}{i}
		(-1)^i
		\sumt[n]{j=i}
		e^{-\frac12 j (\log n + 2\cc)}
		\binomt{n - \abs x}{j - i}
	-	1
	}
+	\oh1
\\&\qquad
=
	2^{-n}
	\sumt{x \in \mbf_2^n}
	\absb{
		\sumt[\abs x]{i=0}
		\binomt{\abs x}{i}
		e^{-\cc i} (-1)^i n^{-i/2}
		\sumt[n]{j=i}
		e^{-\cc(j-i)} n^{-(j-i)/2}
		\binomt{n - \abs x}{j - i}
	-	1
	}
+	\oh1
\\&\qquad
=
	2^{-n}
	\sumt{x \in \mbf_2^n}
	\absb{
		\rbb{ 1 - e^{-\cc} n^{-1/2} }^{\abs x}
		\rbb{ 1 + e^{-\cc} n^{-1/2} }^{n-\abs x}
	-	1
	}
+	\oh1.
\\&\qquad
=
	\sumt[n]{r=0}
	\binomt nr 2^{-n}
	\absb{
		\rbb{ 1 - e^{-\cc/2} n^{-1/2} }^{r}
		\rbb{ 1 + e^{-\cc/2} n^{-1/2} }^{n-r}
	-	1
	}
+	\oh1
\\&\qquad
=
	2 \cdot d_\TV\rbb{ \Bin(n, \tfrac12 - \tfrac12 e^{-\cc} / \sqrt n), \: \Bin(n, \tfrac12) }
+	\oh1.
\]
Applying \cref{res:app:binom-tv} with $\alpha \cq 1$ and $z \cq e^{-\cc}$, we deduce the theorem.
\end{Proof}

\subsection{Total Variation Distance Between Binomials}

In this section of the appendix, we determine a limiting expression for the TV distance between two particular Binomial distributions.
Namely, we prove the following lemma.

\begin{lem}
\label{res:app:binom-tv}
	Let $(n_N)_\Ninn \in \mbn^\mbn$ and $(\aa_N)_\Ninn \in (0,\infty)^\mbn$.
	Suppose that $\min\bra{ \aa_N n_N, \: n_N / \aa_N } \to \infty$ \asinf N.
	Then,
	for all $\yy \in \mbr$,
	we have
	\[
		d_\TV\rbb{
			\Bin\rbr{\nn_N, \: \tfrac{\aa_N}{\aa_N+1} - \tfrac{\aa_N}{\aa_N+1} \yy / \sqrt{\aa_N \nn_N}}, \:
			\Bin(\nn_N, \: \tfrac{\aa_N}{\aa_N+1})
		}
	\to
		2 \, \Phi\rbb{ \tfrac12 \abs{\yy} } - 1
	\quad
		\asinf N.
	\]
\end{lem}

\begin{rmkt*}
	The technical details behind this proof are non-trivial. The statement itself, however, should not be considered deep.
	Indeed, \textcite[Page~59]{DGM:rw-hypercube} need the same result; they simply state, unjustified, that it follow from the CLT for fixed $\aa \in (0,\infty)$.
\end{rmkt*}

As always, we drop the $N$-subscript during the computations in the proof.

\begin{Proof}[Proof of \cref{res:app:binom-tv}]
Our plan is to approximate the Binomial distributions by a discrete normal distribution, using a local CLT, and then approximate this discrete normal by a continuous normal.

We need to set up some notation.
First we explicitly define the distributions.
\begin{itemize}[noitemsep, topsep = \smallskipamount, label = \bcdot]
	\item 
	Write $b_{n,p}$ for the pdf of the $\Bin(n,p)$ distribution:
	\[
		b_{n,p}(k)
	=
		\binomt nk p^k (1-p)^{n-k}
	\Quad{for}
		k \in \bra{0, 1, ..., n}.
	\]
	
	\item 
	Write $\varphi_{\mu,\sigma^2}$ for the pdf of the $\N(\mu,\sigma^2)$ distribution:
	\[
		\varphi_{\mu,\sigma^2}(x)
	=
		\rbb{ 2 \pi \sigma^2 }^{-1/2} \expb{ - \tfrac1{2\sigma^2} \rbr{ x - \mu }^2 }
	\Quad{for}
		x \in \mbr.
	\]
\end{itemize}
Now we choose the parameters for these distributions.
\begin{itemize}[noitemsep, topsep = \smallskipamount, label = \bcdot]
	\item 
	Set $p \cq \tfrac{\aa}{\aa+1}$;
	set $\tilde p \cq \tfrac{\aa}{\aa+1} - \tfrac{\aa}{\aa+1} \yy / \sqrt{\aa n}$.
	
	\item 
	Set $\mu \cq pn = \tfrac{\aa}{\aa+1} n$ and $\sigma^2 \cq p(1-p)n = \tfrac{\aa}{(\aa+1)^2} n$;
	set $\tilde \mu \cq \tilde p n$ and $\tilde \sigma^2 \cq \tilde p (1 - \tilde p) n$.
\end{itemize}
These parameters are related in the following way:
\[
	\rbr{ \mu - \tilde \mu } / \sigma
&
=
	\rbb{ n \tfrac{\aa}{\aa+1} \yy/\sqrt{\aa n}} \big/ \rbb{ n \tfrac{\aa}{(\aa+1)^2} }^{1/2}
=
	\yy;
\\
	\tilde \sigma^2 / \sigma^2
&
=
	1 + \yy (1 - \tfrac1\aa) \sqrt{\aa/n} - \yy^2/n.
\]
In order to apply a local CLT, we need to restrict the distributions to an interval on which asymptotically all the mass is supported: set
\[
	I_n
\cq
	\sbb{
		\tfrac{\aa}{\aa+1} n - \omega \rbr{ \tfrac{\aa}{\aa+1} n }^{1/2}, \:
		\tfrac{\aa}{\aa+1} n + \omega \rbr{ \tfrac{\aa}{\aa+1} n }^{1/2}
	} \cap \mbz
\]
where $\omega \gg 1$ diverges arbitrarily slowly.
Also write $d'_\TV$ to indicate TV distance between two distributions, but restricted to $I_n$.
Now, if $\nu$ is any of the above distributions, then $\nu(I_n) = 1 - \oh1$.
Hence, for any two such distributions $\nu$ and $\tilde \nu$, we have
\(
	d_\TV(\nu,\tilde \nu) = d'_\TV(\nu,\tilde \nu) + \oh1.
\)

Before calculating the TV, we make some preliminary approximations.
It is easy to check that
\[
	\MAX{x \in I_n}
	\MAX{\xi \in [-\frac12,\frac12]}
	\absb{ \varphi_{\mu,\sigma^2}(x + \xi) - \varphi_{\mu,\sigma^2}(x) } / \varphi_{\mu,\sigma^2}(x)
=
	\oh1.
\]
It is also easy to check, using Stirling's approximation, that
\[
	\delta_n(k)
\cq
	b_{n,p}(k) / \varphi_{\mu,\sigma^2}(k) - 1
\Quad{for}
	k \in \mbz
\Quad{satisfies}
	\LIM{\toinf n}
	\MAX{k \in I_n}
	\: \abs{ \delta_n(k) }
=
	0.
\]
Analogous results hold when $(p,\mu,\sigma^2)$ is replaced by $(\tilde p, \tilde \mu, \tilde \sigma^2)$, defining $\tilde \delta_n(\cdot)$ similarly.

Having done all this preparation, we are eventually ready to calculate the TV in question:
\[
&	d_\TV\rbb{ \Bin(n, p), \: \Bin(n, \tilde p) }
=
	d'_\TV\rbb{ \Bin(n, p), \: \Bin(n, \tilde p) }
+ \oh1
\\&\qquad
=
	\tfrac12
	\sumt{k \in I_n}
	\absb{ b_{n,p}(k) - b_{n,\tilde p}(k) }
=
	\tfrac12
	\sumt{k \in I_n}
	\absb{ \varphi_{\mu,\sigma^2}(k) - \varphi_{\tilde \mu, \tilde \sigma^2}(k) }
+	\oh1
\\&\qquad
=
	\tfrac12
	\intt{I_n}
	\absb{ \varphi_{\mu,\sigma^2}(x) - \varphi_{\tilde \mu, \tilde \sigma^2}(x) } \, dx
+	\oh1
=
	d_\TV\rbb{ \N(\mu, \sigma^2), \: \N(\tilde \mu, \tilde \sigma^2) }
+	\oh1.
\]
It remains to calculate this TV distance between two normal distributions, for which we have a nice pdf.
First,
by translation and scaling,
it is straightforward to see that
\[
	d_\TV\rbb{ \N(\mu, \sigma^2), \: \N(\tilde \mu, \tilde \sigma^2) }
=
	d_\TV\rbb{ \N(0, 1), \: \N( \yy, 1 + \yy (1 - \tfrac1\aa) \sqrt{\aa/n} - \yy^2 / n) }.
\]
Next, we claim that
\[
	d_\TV\rbb{ \N(0, 1), \: \N(0, 1 + \eps) } \to 0
\quad
	\aszero \eps.
\]
Applying this, translated by $\yy$,
with $\eps \cq \yy (1 - \tfrac1\aa) \sqrt{\aa/n} - \yy^2/n$,
we obtain
\[
	d_\TV\rbb{ \N(\mu, \sigma^2), \: \N(\tilde \mu, \tilde \sigma^2) }
=
	d_\TV\rbb{ \N(0,1), \: \N(\yy,1) }
+	\oh1
\Quad{provided}
	1/n \ll \aa \ll n.
\]
Manipulating integrals, using the fact that $\varphi_{m,1}(x) = \varphi_{0,1}(x-m)$ for any $x,m \in \mbr$,
writing
\[
	\Phi_{0,1}(\beta) \cq \intt[\beta]{-\infty} \varphi_{0,1}(x) \, dx
\Quad{for}
	\beta \in \mbr
\]
for the standard normal cumulative density function, we find that
\[
	d_\TV\rbb{ \N(0,1), \: \N(\yy,1) }
=
	2 \, \Phi_{0,1}\rbb{\tfrac12 \abs{\yy}} - 1.
\]
(In fact, this is not specialised to the normal distribution: it works for many distributions.)
Thus
\[
	d_\TV\rbb{ \Bin(n, p), \: \Bin(n, \tilde p) }
\to
	2 \, \Phi_{0,1}\rbb{\tfrac12 \abs{\yy}} - 1
\quad
	\asinf n.
\]

It remains to prove that
\[
	d_\TV\rbb{ \N(0, 1), \: \N(0, 1 - \eps) } \to 0
\quad
	\aszero \eps.
\]
To see this, first observe that,
for $\eps \in [-\tfrac1{10}, \tfrac1{10}]$ and $x \in \mbr$,
we have
\[
	\varphi_{0,1/(1+\eps^5)}(x) / \varphi_{0,1}(x)
=
	\expb{ - \tfrac12 \eps^5 x^2}.
\]
We now split the integral
\(
	\intt[\infty]{-\infty}
=
	\intt[-1/\eps]{-\infty} + \intt[1/\eps]{-1/\eps} + \intt[\infty]{1/\eps}:
\)
in the middle region, we have
\[
	\abs{ \varphi_{0,1/(1+\eps^4)}(x) / \varphi_{0,1}(x) - 1 }
\le
	\abs{ \expb{ \tfrac12 \abs \eps^3 } - 1 }
\le
	\eps^3;
\]
the probability that either random variable lands in the outer regions tends to 0 as $\eps \to 0$.
Hence
\[
	d_\TV\rbb{ N(0,1), \: N(0,1/(1+\eps^5)) }
\to
	0
\quad
	\aszero \eps.
\]
This proves the stated claim, and hence completes the proof of the lemma.
\end{Proof}

\end{document}